\documentclass[aap,MSNbibl,seceqn,citesort,dvips]{arximspdf}

\doi{10.1214/11-AAP817} 
\volume{22}
\issue{5}
\pubyear{2012}
\firstpage{1880}
\lastpage{1927}

\makeatletter

\newtheorem{theorem}{Theorem}[section]
\newtheorem{coro}[theorem]{Corollary}
\newtheorem{lem}[theorem]{Lemma}

\newcommand{\cal}{\mathcal}

\newcommand{\eqd}{\stackrel{D}{=}}
\newcommand{\convd}{\stackrel{D}{\longrightarrow}}
\newcommand{\convp}{\stackrel{p}{\longrightarrow}}
\newcommand{\dr}{d \rightarrow\infty}
\newcommand{\nr}{n \rightarrow\infty}
\newcommand{\ez}{\mathbb{E}}
\newcommand{\pz}{\mathbb{P}}

\newcommand{\xx}{\hat{\mathbf{X}}}

\makeatother

\begin{document}
\begin{frontmatter}

\title{Optimal scaling of random walk Metropolis algorithms with discontinuous target densities}
\runtitle{Optimal scaling of RWM algorithms}

\begin{aug}
\author[A]{\fnms{Peter} \snm{Neal}\corref{}\ead[label=e1]{P.Neal-2@manchester.ac.uk}},
\author[B]{\fnms{Gareth} \snm{Roberts}\ead[label=e2]{Gareth.O.Roberts@warwick.ac.uk}}
\and
\author[C]{\fnms{Wai Kong} \snm{Yuen}\ead[label=e3]{jyuen@brocku.ca}}
\runauthor{P. Neal, G. Roberts and W. K. Yuen}
\affiliation{University of Manchester, University of Warwick and
Brock University}
\address[A]{P. Neal\\
School of Mathematics \\
University of Manchester \\
Sackville Street \\
Manchester\\
M60 1QD \\
United Kingdom \\
\printead{e1}}
\address[B]{G. Roberts\\
Department of Statistics \\
University of Warwick \\
Coventry \\
CV4 7AL \\
United Kingdom \\
\printead{e2}}
\address[C]{W. K. Yuen\\
Department of Mathematics \\
Brock University \\
St. Catharines \\
ON L2S 3A1 \\
Canada \\
\printead{e3}} 
\end{aug}

\received{\smonth{2} \syear{2011}}
\revised{\smonth{9} \syear{2011}}

%
\begin{abstract}
We consider the optimal scaling problem for high-dimensional random
walk Metropolis (RWM) algorithms where the target distribution has a
discontinuous probability density function. Almost all previous
analysis has focused upon continuous target densities. The main
result is a weak convergence result as the dimensionality $d$ of the
target densities converges to $\infty$. In particular, when the
proposal variance is scaled by $d^{-2}$, the sequence of stochastic
processes formed by the first component of each Markov chain
converges to an appropriate Langevin diffusion process. Therefore
optimizing the efficiency of the RWM algorithm is equivalent to
maximizing the speed of the limiting diffusion. This leads to an
asymptotic optimal acceptance rate of $e^{-2}$ $(\mbox{$=$}0.1353)$ under quite
general conditions. The results have major practical implications
for the implementation of RWM algorithms by highlighting the
detrimental effect of choosing RWM algorithms over
Metropolis-within-Gibbs algorithms.
\end{abstract}

%

\end{frontmatter}

\section{Introduction} \label{SecInt}

Random walk Metropolis (RWM) algorithms are widely used generic
Markov chain Monte Carlo (MCMC) algorithms. The ease with which RWM
algorithms can be constructed has no doubt played a pivotal role in
their popularity. The efficiency of a RWM algorithm depends
fundamentally upon the scaling of the proposal density. Choose the
variance of the proposal to be too small and the RWM will converge
slowly since all its increments are small. Conversely, choose the
variance of the proposal to be too large and too high a proportion
of proposed moves will be rejected. Of particular interest is how
the scaling of the proposal variance depends upon the dimensionality
of the target distribution. The target distribution is the
distribution of interest and the MCMC algorithm is constructed such
that the stationary distribution of the Markov chain is the target
distribution.

The Introduction is structured as follows. We outline known results for
continuous independent and identically distributed product densities
from \cite{RGG} and subsequent work. We highlight the scope and
limitations of the results before introducing the discontinuous target
densities to be studied in this paper. While the statements of the key
results (Theorem~\ref{main}) in this paper are similar to those given
for continuous target densities, the proofs are markedly different. A
discussion of why a new method of proof is required for discontinuous
target densities is given. Finally, we give an outline of the
remainder of the paper.

The results of this paper have quite general consequences for the
implementation of Metropolis algorithms on discontinuous densities (as
are commonly applied in many Bayesian Statistics problems), namely:
\begin{longlist}[(1)]
\item[(1)] Full- (high-) dimensional update rules can be an order of
magnitude slower than strategies involving smaller dimensional updates.
(See Theorem~\ref{thmprop} below.)
\item[(2)] For target densities with bounded support, Metropolis algorithms
can be an order of magnitude slower than algorithms which first
transform the target support to $\mathbb{R}^d$ for some $d$.
\end{longlist}

In~\cite{RGG}, a sequence of target densities of the form
%
%
\begin{equation} \label{eqrevb1} \pi_d (\mathbf{x}^d) = \prod
_{i=1}^d f (x_i^d)
\end{equation}
were considered as $\dr$, where $f(\cdot)$ is twice differentiable
and satisfies certain mild moment conditions; see~\cite{RGG}, (A1)
and (A2). The following random walk Metropolis algorithm was used to
obtain a sample $\mathbf{X}_0^d, \mathbf{X}_1^d, \ldots$ from $\pi_d
(\cdot)$. Draw $\mathbf{X}_0^d$ from $\pi_d (\cdot)$. For $t \geq0$
and $i=1,2, \ldots,$ let $Z_{t,i}$ be independent and identically
distributed (i.i.d.) according to $Z \sim N(0,1)$ and
$\mathbf{Z}_t^d = (Z_{t,1}, Z_{t,2}, \ldots, Z_{t,d})$. At time $t$,
propose
%
%
\begin{equation} \label{eqrevb2} \mathbf{Y}^d = \mathbf{X}_t^d +
\sigma_d \mathbf{Z}_t^d,
\end{equation}
where $\sigma_d$ is the proposal standard deviation to be discussed
shortly. Set $\mathbf{X}_{t+1}^d = \mathbf{Y}^d$ with probability
%
%
\begin{equation} \label{eqrevb3} \alpha(\mathbf{X}_t^d, \mathbf
{Y}^d) \equiv1 \wedge
\frac{\pi_d (\mathbf{Y}^d)}{\pi_d (\mathbf{X}_t^d)}.
\end{equation}
Otherwise set $\mathbf{X}_{t+1}^d = \mathbf{X}_t^d$. It is
straightforward to check that $\{ \mathbf{X}_t^d \}$ has stationary
distribution $\pi_d (\cdot)$, and hence, for all $t \geq0$,
$\mathbf{X}_t^d \sim\pi_d (\cdot)$. The key question addressed in
\cite{RGG} was: starting from the stationary distribution, how
should $\sigma_d$ be chosen to optimize the rate at which the RWM
algorithm explores the stationary distribution? Since the components
of $\mathbf{X}_t^d$ are i.i.d., it suffices to study the marginal
behavior of the first component, $X_{t,1}^d$. In~\cite{RGG}, it was
shown that if $\sigma_d = l/\sqrt{d}$ $(l > 0)$ and $U_t^d= X_{[t d],1}^d$
$(t \geq0)$, then
%
%
\begin{equation} \label{eqrevb4} U^d \Rightarrow U
\qquad\mbox{as } \dr,
\end{equation}
where $U_\cdot$ satisfies the Langevin SDE
%
%
\begin{equation} \label{eqrevba5} dU_t = \sqrt{h (l)} \,dB_t + \phi
(l) \frac{f^\prime(U_t)}{2 f (U_t)} \,dt
\end{equation}
with $U_0 \sim f (\cdot)$ and $h(l) = 2 l^2 \Phi(-l \sqrt{I}/2)$
with $\Phi$ being the standard normal c.d.f. and $I \equiv\ez_f [ \{
f^\prime(X)/f (X) \}^2]$. Note that the ``speed measure'' of the
diffusion $\phi(l)$ only depends upon $f$ through $I$. The
diffusion limit for $U^d$ is unsurprising in that for a time interval of
length $s>0$, $O(d)$ moves are made each of size $O(1/\sqrt{d})$.
Therefore the movements in the first component (appropriately
normalized) converge to those of a Langevin diffusion with the ``most
efficient'' asymptotic diffusion having the largest speed measure~$h(l)$.
Since the diffusion limit involves speeding up time by a
factor of~$d$, we say that the mixing of the algorithm is $O(d)$.
The optimal value of~$l$ is $\hat{l} = 2.38/\sqrt{I}$, which leads
to an average optimal acceptance rate (aoar) of 0.234. This has
major practical implications for practitioners, in that, to monitor
the (asymptotic) efficiency of the RWM algorithm it is sufficient to
study the proportion of proposed moves accepted.

There are three key assumptions made in~\cite{RGG}. First,
$\mathbf{X}_0^d \sim\pi_d (\cdot)$, that is, the algorithm starts
in the stationary distribution and $\sigma_d$ is chosen to optimize
exploration of the stationary distribution. This assumption has been
made in virtually all subsequent optimal scaling work; see, for
example,~\cite{BR00,NR08,NR06,Bed06} and
\cite{RR01}. The one exception is~\cite{CRR05}, where
$\mathbf{X}_0^d$ is started from the mode of $\pi_d (\cdot)$ with
explicit calculations given for a standard multivariate normal
distribution. In~\cite{CRR05}, it is shown that $\sigma_d =
O(1/\sqrt{d})$ is optimal for maximizing the rate of convergence to
the stationary distribution. Since convergence is shown to occur
within $O( \log d)$ iterations, the time taken to explore the
stationary distribution dominates the time taken to converge to the
stationary distribution, and thus overall it is optimal to choose
$\sigma_d = \hat{l}/\sqrt{d}$. It is difficult to prove generic
results for $\mathbf{X}_0^d \not\sim\pi_d$. However, the findings
of~\cite{CRR05} suggest that even when $\mathbf{X}_0^d \not\sim
\pi_d$, it is best to scale the proposal distribution based upon
$\mathbf{X}_0^d \sim\pi_d$. It is worth noting that in~\cite{CRR05}
it was found that for the Metropolis adjusted Langevin algorithm
(MALA), the optimal scaling of $\sigma_d$ for $\mathbf{X}_0^d$
started at the mode of a multivariate normal is $O(d^{-1/4})$
compared to $O(d^{-1/6})$ for $\mathbf{X}_0^d \sim\pi_d$.

Second, $\pi_d (\cdot)$ is an i.i.d. product density. This
assumption has been relaxed by a number of authors with $\sigma_d =
O(1/\sqrt{d})$ and an aoar of 0.234 still being the case, for
example, independent, scaled product densities (\cite{RR01}
and~\cite{Bed06}), Gibbs random fields~\cite{BR00}, exchangeable
normals~\cite{NR06} and elliptical densities~\cite{SR09}. Thus
the simple rule of thumb of tuning $\sigma_d$ such that one in four
proposed moves are accepted holds quite generally. In~\cite{Bed08}
and~\cite{SR09}, examples where the aoar is strictly less than 0.234
are given. These correspond to different orders of magnitude being
appropriate for the scaling of the proposed moves in different
components.

Third, the results are asymptotic as $\dr$. However, simulations
have shown that for i.i.d. product densities an acceptance rate
of 0.234 is close to optimal for $d=10$; see, for example,
\cite{NR06}. Departures from the i.i.d. product density require
larger $d$ for the asymptotic results to be optimal, but $d=100$ is
often seen in practical MCMC problems. In~\cite{NR10} and
\cite{Sherlock}, optimal acceptance rates are obtained for finite $d$
for some special cases.

With the exceptions of~\cite{NR08,NR10} and~\cite{SR09}, in
the above works $\pi_d$ is assumed to have a continuous (and
suitably differentiable) probability density function (p.d.f.). The aim
of the current work is to investigate the situation where the target
distribution has a discontinuous p.d.f., and specifically, target
distributions confined to the $d$-dimensional hypercube $[0,1]^d$.
That is, we consider target distributions of the form
%
%
\begin{equation} \label{eq1a} \pi_d (\mathbf{x}^d) = \prod_{i=1}^d f
(x_i^d),
\end{equation}
where
%
%
\begin{equation} \label{eq1b} f (x) \propto
\exp(g(x)) 1_{\{0 < x < 1\}} \qquad(x \in\mathbb{R})
\end{equation}
and $g ( \cdot)$
is twice differentiable upon $[0,1]$ with
%
%
\begin{equation} \label{eq1c} g^\ast= {\sup_{0 \leq y
\leq1}} |g^\prime(y)| < \infty.
\end{equation}
We then use the following random
walk Metropolis algorithm to obtain a~sample $\mathbf{X}_0^d,
\mathbf{X}_1^d, \ldots$ from $\pi_d (\cdot)$. Draw $\mathbf{X}_0^d$
from $\pi_d (\cdot)$. For $t \geq0$ and $i=1,2, \ldots,$ let
$Z_{ti}$ be independent and identically distributed (i.i.d.)
according to $Z \sim U[-1,1]$ and $\mathbf{Z}_t^d = (Z_{t1}, Z_{t2},
\ldots, Z_{td})$. At time $t$, propose
%
%
\begin{equation} \label{eqrevb5} \mathbf{Y}^d = \mathbf{X}_t^d +
\sigma_d \mathbf{Z}_t^d.
\end{equation}
Set $\mathbf{X}_{t+1}^d = \mathbf{Y}^d$ with probability
%
%
\begin{equation} \label{eqrevb6} \alpha(\mathbf{X}_t^d, \mathbf
{Y}^d) \equiv1 \wedge
\frac{\pi_d (\mathbf{Y}^d)}{\pi_d (\mathbf{X}_t^d)}.
\end{equation}
Otherwise set $\mathbf{X}_{t+1}^d = \mathbf{X}_t^d$.

In~\cite{NR08} and~\cite{SR09}, spherical and elliptical densities
are considered which have very different geometry to the hypercube
restricted densities. Therefore different approaches are taken in
these papers with results akin to those obtained for continuous
target densities. Densities of the form (\ref{eq1b}) have previously
been studied in~\cite{NR10}, where the expected square jumping
distance (ESJD) has been computed. The ESJD is
%
%
\begin{equation} \label{eqreve1} \ez_{\pi_d} \biggl[ \sum_{i=1}^d
(X_{1,i}^d - X_{0,i}^d)^2
\biggr] = d \ez_{\pi_d}[ (X_{1,1}^d - X_{0,1}^d)^2 ],
\end{equation}
the mean squared distance between $\mathbf{X}_0^d$
and $\mathbf{X}_1^d$, where $\mathbf{X}_0^d \sim\pi_d$. In
\cite{NR10}, Appendix~B, it is shown that for $\sigma_d = l/d$ $(l >
0)$ and $f (x) = 1_{\{0 < x < 1 \}}$,
%
%
\begin{equation} \label{eqreve2} d \ez_{\pi_d} \Biggl[ \sum
_{i=1}^d (X_{1,i}^d - X_{0,i}^d)^2
\Biggr] \rightarrow\frac{l^2}{3} \exp\biggl( - \frac{l}{2} \biggr)
\qquad\mbox{as } \dr.
\end{equation}
Thus asymptotically
(as $\dr$) the ESJD is maximized by taking $\hat{l}=4$ which
corresponds to an aoar of $\exp(-2)$ ($=$0.1353). In this paper, we
show that $\sigma_d = l/d$ and an aoar of $\exp(-2)$ holds more
generally for target distributions of the form given by (\ref{eq1a})
and (\ref{eq1b}). Moreover, we prove a much stronger result than
that given in~\cite{NR10}, in that, we prove that $V^d_s = X_{[s
d^2],1}^d$ converges weakly to an appropriate Langevin diffusion
$V_s$ with speed measure $\phi(l) =(l^2/3) \exp(-l/(2f^\ast))$ as
$\dr$, where $f^\ast= \lim_{x \downarrow0} \{ (f(x) + f(1-x))/2
\}$. This gives a clear indication of how the Markov chain explores
the stationary distribution. By contrast the ESJD only gives a
measure of average behavior and does not take account of the
possibility of the Markov chain becoming ``stuck.'' If
$\ez_{\mathbf{Z}^d} [ \alpha(\mathbf{x}^d, \mathbf{x}^d+
\mathbf{Z}^d)]$ is very low, the Markov chain started
$\mathbf{X}_0^d = \mathbf{x}^d$ is likely to spend a large number of
iterations at $\mathbf{x}^d$ before accepting a move away from
$\mathbf{x}^d$. Note that since $V_s^d$ involves speeding up time by
a factor of $d^2$, we say that the mixing of the algorithm is
$O(d^2)$. The ESJD is easy to compute and asymptotically, as $\dr$,
the ESJD (appropriately scaled) converges to $\phi(l)$. Thus in
discussing possible extensions of the Langevin diffusion limit
proved in Theorem~\ref{main} for i.i.d. product densities of the
form given in (\ref{eq1a}) and (\ref{eq1b}), we make considerable
use of the ESJD. However, we highlight the limitations of the ESJD
in discussing extensions of Theorem~\ref{main}.

In most previous work on optimal scaling, the components of
$\mathbf{Z}^d$ are taken to be independent and identically
distributed $Z \sim N(0,1)$ random variables. The reason for
choosing $Z \sim U[-1,1]$ for discontinuous target densities is
mathematical convenience. The results proved in this paper hold with
Gaussian rather than uniform proposal distributions, but some
elements of the proof are less straightforward. For discussion of
the ESJD for densities (\ref{eq1a}) for general $Z$ subject to $\ez
[Z^2] <\infty$, see~\cite{NR10}, Appendix B.

While the key result, a Langevin diffusion limit for the movement
in the first component, is the same as~\cite{RGG}, the proof is
markedly different. Note that, for finite $d$, $U^d$ and $V^d$ are
not Markov chains since whether or not a proposed move is accepted
depends upon all the components in $\pi_d (\cdot)$. In~\cite{RGG},
it is shown that there exists $\{ F_d \}$ such that $\pz(
\bigcup_{t=0}^{[Td]} \{\mathbf{X}_t^d \notin F_d \}) \rightarrow0$ as
$\dr$ and
%
%
\begin{equation} \label{eqrevb7} \sup_{\mathbf{x}^d \in F_d}
\biggl| \ez
[ \alpha(\mathbf{x}^d, \mathbf{x}^d + \sigma_d \mathbf{Z}^d)] -
2\Phi\biggl( -\frac{l \sqrt{I}}{2} \biggr) \biggr| \leq\varepsilon_d,
\end{equation}
where $\varepsilon_d \rightarrow0$ as $\dr$. While (\ref{eqrevb7})
is not explicitly stated in~\cite{RGG}, it is the essence of the
requirements of the sets $\{F_d\}$, stating that for large $d$, with
high probability over the first $Td$ iterations the acceptance
probability of the Markov chain is approximately constant, being
within $\varepsilon_d$ of $2 \Phi(-l \sqrt{I}/2)$.
(Note $n$ rather than $d$ is used for dimensionality in~\cite{RGG}.)
Thus in the limit as $\dr$ the effect of the other components on
movements in the first component converges to a deterministic
acceptance probability $2 \Phi(-l \sqrt{I}/2)$. The situation is
more complex for $\pi_d (\cdot)$ of the form given by (\ref{eq1a})
and (\ref{eq1b}) as the acceptance rate in the limit as $\dr$ is
inherently stochastic. For example, suppose $\pi_d (\cdot)$ is the
uniform distribution on the $d$-dimensional hypercube so that
$\alpha(\mathbf{X}_t^d, \mathbf{Y}^d) = 1_{\{ \mathbf{Y}^d \in
[0,1]^d \}}$. Letting $R_d^L = (0, \sigma_d)$ and $R_d^U =
(1-\sigma_d,1)$, this gives
%
%
\begin{equation} \label{eqrevb8}
\ez[ \alpha(\mathbf{x}^d, \mathbf{x}^d + \sigma_d \mathbf{Z}^d)] =
\prod_{i \in R_d^L} \biggl( \frac{1}{2} + \frac{x_i}{2 \sigma_d}
\biggr) \times\prod_{i \in R_d^U} \biggl( \frac{1}{2} +
\frac{1-x_i}{2 \sigma_d} \biggr).
\end{equation}
Thus the acceptance probability is totally determined by the
components at the boundary (within $\sigma_d$ of 0 or 1). The total
number of components in $R_d^L \cup R_d^U$ is $\operatorname{Bin} (d,2l/d)$
which converges in distribution to $\operatorname{Po} (2l)$ as $\dr$. Thus
the number of components close to the boundary is inherently
stochastic. Moreover, the location of the components within $R_d^L
\cup R_d^U$ plays a crucial role in the acceptance probability; see
(\ref{eqrevb8}). Therefore there is no hope of replicating
directly the method of proof applied in~\cite{RGG} and subsequently,
in~\cite{BR00} and~\cite{NR06}.

We need a homogenization argument which involves looking at
$\mathbf{X}^d_\cdot$ over $[d^\delta]$ steps; cf.~\cite{NR08}. In
particular, we show that the acceptance
probability converges very rapidly to its stationary measure, so
that over $[d^\delta]$ iterations approximately $\exp(-lf^\ast/2)
[d^\delta]$ proposed moves are accepted. By comparison, $|
X_{[d^\delta],1}^d - X_{0,1}^d | \leq[d^\delta] \sigma_d$; thus the
value of an individual component only makes small changes over
$[d^\delta]$ iterations. That is, we show that there exists $\{
\tilde{F}_d \}$ such that, for any $T >0$, $\pz(
\bigcup_{t=0}^{[Td^2]} \{\mathbf{X}_t^d \notin\tilde{F}_d \})
\rightarrow0$ as $\dr$ and for $\delta>0$,
%
%
\begin{equation} \label{eqrevb10} \sup_{\mathbf{x}^d \in\tilde
{F}_d} \Biggl| \frac{1}{[d^\delta]} \sum_{t=0}^{
[d^\delta]-1} \ez[ \alpha(\mathbf{X}^d_t, \mathbf{X}^d_t +
\sigma_d \mathbf{Z}_t^d) | \mathbf{X}_0^d = \mathbf{x}^d] - \exp
\biggl( - \frac{lf^\ast}{2} \biggr) \Biggr| \leq\varepsilon_d\hspace*{-25pt}
\end{equation}
for some $\varepsilon_d \rightarrow0$ as $\dr$. For
large $d$, with high probability over the first $[Td^2]$ iterations
the Markov chain stays in $\tilde{F}_d$, where the average number of
accepted proposed moves in the following $[d^\delta]$ iterations is
$\exp(-lf^\ast/2)d^\delta+ o(d^\delta)$. The arguments are
considerably more involved than in~\cite{NR08}, where spherically
constrained target distributions were studied, due to the very
different geometry of the hypercube and spherical constraints
applied in this paper and~\cite{NR08}, respectively. In particular,
in~\cite{NR08}, $\sigma_d = l/\sqrt{d}$ with an aoar of $0.234$.

By exploiting the homogenization argument it is possible to prove
that~$V^d$ converges weakly to an appropriate Langevin diffusion
$V$, given in Theorem~\ref{main}. In Section~\ref{SecAlg}, Theorem
\ref{main} is presented along with an outline of the proof. Also in
Section~\ref{SecAlg}, a description of the pseudo-RWM algorithm is
given. The pseudo-RWM algorithm plays a key role in the proof of
Theorem~\ref{main}. The pseudo-RWM process moves at each iteration
and the moves in the pseudo-RWM process are identical to those of
the RWM process, conditioned upon a proposed move in the RWM process
being accepted. The proof of Theorem~\ref{main} is long and
technical with the details split into three key sections which are
given in the \hyperref[app]{Appendix}; see Section~\ref{SecAlg} for
more details. In
Section~\ref{SecExtSim}, two interesting extensions of Theorem
\ref{main} are given. In particular, Theorem~\ref{thmprop} has major
practical implications for the implementation of RWM algorithms by
highlighting the detrimental effect of choosing RWM algorithms over
Metropolis-within-Gibbs algorithms. The target densities for which
theoretical results can be proved are limited, so discussion of
possible extensions of Theorem~\ref{main} are given. In particular,
we discuss general~$\pi_d$ restricted to the hypercube, general
discontinuities in $f$ and $\mathbf{X}_0^d \not\sim\pi_d$.

\section{\texorpdfstring{Pseudo-RWM algorithm and Theorem \protect\ref{main}}{Pseudo-RWM algorithm and Theorem 2.1}} \label{SecAlg}

We begin by defining the pseudo-random walk Metropolis (pseudo-RWM)
process. We will then be in position to formally state the main
theorem, Theorem~\ref{main}. An outline of the proof of Theorem
\ref{main} is given, with the details, which are long and technical,
placed in the \hyperref[app]{Appendix}.

For $d \geq1$, let
\[
h_d (\mathbf{z}^d) = \cases{
2^{-d}, &\quad if $\mathbf{z}^d \in(-1,1)^d$, \cr
0, &\quad otherwise.}
\]
Let $J_d (\mathbf{x}^d)$ denote the probability of
accepting a move in the RWM process given the current state of the
process is $\mathbf{x}^d$. Then
%
%
\begin{equation} \label{eq21} J_d
(\mathbf{x}^d) = \int h_d (\mathbf{z}^d ) \biggl\{ 1 \wedge
\frac{\pi_d (\mathbf{x}^d + \sigma_d \mathbf{z}^d)}{\pi_d
(\mathbf{x}^d)} \biggr\} \,d \mathbf{z}^d.
\end{equation}
Let\vspace*{1pt} $b_d^r (\mathbf{x}^d) = \sum_{j=1}^d 1_{ \{ x_j \in R_d^r \}}$,
the total number of components of $\mathbf{x}^d$ in $R_d^r =(0,r/d)
\cup(1-r/d,1)$. By Taylor's theorem for all $0 \leq x_i, x_i +
\sigma_d z_i \leq1$ and $-1 \leq z_i \leq1$,
%
%
\begin{equation} \label{eq21az} g(x_i + \sigma_d z_i) - g(x_i) \geq
-g^\ast\sigma_d
\end{equation}
with $g^\ast$ defined in (\ref{eq1c}).
Hence, for all $\mathbf{x}^d \in[0,1]^d$,
%
%
\begin{eqnarray} \label{eq21ay} J_d
(\mathbf{x}^d) & = & \int h_d (\mathbf{z}^d) \Biggl\{ 1 \wedge
\prod_{i=1}^d \frac{\exp( g (x_i + \sigma_d z_i))}{
\exp(g(x_i))} \Biggr\} 1_{\{ \mathbf{x}^d + \sigma_d \mathbf
{z}^d \in[0,1]^d \}} \,d \mathbf{z}^d \nonumber\\
& \geq&\int h_d (\mathbf{z}^d) \{ 1 \wedge\exp(- dg^\ast
\sigma_d) \} 1_{\{ \mathbf{x}^d + \sigma_d \mathbf
{z}^d \in
[0,1]^d \}} \,d \mathbf{z}^d \nonumber\\[-8pt]\\[-8pt]
&= & \exp(-l g^\ast) \int h_d (\mathbf{z}^d) 1_{\{ \mathbf
{x}^d +
\sigma_d \mathbf{z}^d \in[0,1]^d \}} \,d \mathbf{z}^d
\nonumber\\
& \geq& \exp(-l g^\ast) \biggl( \frac{1}{2} \biggr)^{b_d^l
(\mathbf{x}^d)}.\nonumber
\end{eqnarray}
This lower bound for $J_d
(\mathbf{x}^d)$ will be used repeatedly.

The pseudo-RWM process moves at each iteration, which is the key
difference to the RWM process. Furthermore, the moves in the
pseudo-RWM process are identical to those of the RWM process,
conditioned upon a move in the RWM process being accepted, that is, its
jump chain. For $d
\geq1$, let $\hat{\mathbf{X}}_0^d, \hat{\mathbf{X}}_1^d, \ldots$
denote the successive states of the pseudo-RWM process, where
$\hat{\mathbf{X}}_0^d \sim\pi_d (\cdot)$. The pseudo-RWM process is
a Markov process, where for $t \geq0$, $\hat{\mathbf{X}}_{t+1}^d =
\hat{\mathbf{X}}_t^d + \sigma_d \hat{\mathbf{Z}}_t^d$ and given that
$\hat{\mathbf{X}}_t^d = \mathbf{x}^d$, $\hat{\mathbf{Z}}_t^d$ has
p.d.f.
\[
\zeta( \mathbf{z}^d | \mathbf{x}^d) = h_d (\mathbf{z}^d)
\alpha(\mathbf{x}^d, \mathbf{x}^d + \sigma_d \mathbf{z}^d) /
J_d (\mathbf{x}^d), \qquad\mathbf{z}^d \in(-1,1)^d.
\]
Note that $\zeta( \mathbf{z}^d | \mathbf{x}^d) =0$ for
$\mathbf{z}^d \notin(-1,1)^d$. Since $\mathbf{X}_0^d,
\hat{\mathbf{X}}_0^d \sim\pi_d$, we can couple the two processes to
have the same starting value $\mathbf{X}_0^d$. A continued coupling
of the two processes is outlined below. Suppose that $\mathbf{X}_t^d
= \mathbf{x}^d$. Then for any $s \geq1$,
%
%
\begin{equation}
\label{eqn111} \pz\Biggl( \bigcup_{j=1}^s \{ \mathbf{X}_{t+j}^d =
\mathbf{x}^d \} | \mathbf{X}_t^d = \mathbf{x}^d \Biggr) = \bigl(1 - J_d
(\mathbf{x}^d)\bigr)^s.
\end{equation}
That is, the number of iterations the RWM algorithm stays at
$\mathbf{x}^d$ before moving follows a geometric distribution with
``success'' probability $J_d (\mathbf{x}^d)$. Therefore for $j \geq
0$, let
$M_j (\cdot)$ denote independent geometric random
variables, where for $0 < p \leq1$, $M_j (p)$ denotes a geometric
random variable with ``success'' probability $p$. For $s \in
\mathbb{Z}^+$, let $\hat{M}_s^d = M_s (J (\hat{\mathbf{X}}_s^d))$
and for $t \in\mathbb{Z}^+$, let
\[
U_t^d = \sup\Biggl\{ s \in\mathbb{Z}^+\dvtx\sum_{j=0}^{s-1} M_j (J_d
(\hat{\mathbf{X}}_j^d)) \leq t \Biggr\},
\]
where the sum is zero if vacuous.
For $s \in\mathbb{Z}^+$, attach $\hat{M}_s^d = M_s (J
(\hat{\mathbf{X}}_s^d))$ to~$\hat{\mathbf{X}}_s^d$. Thus
$\hat{M}_s^d$ denotes the total number of iterations the RWM process
spends at $\hat{\mathbf{X}}_s^d$ before moving to\vadjust{\goodbreak}
$\hat{\mathbf{X}}_{s+1}^d$. Hence, the RWM process can be
constructed from $(\hat{\mathbf{X}}_0^d, \hat{M}_0^d),
(\hat{\mathbf{X}}_1^d, \hat{M}_1^d), \ldots$ by setting
$\mathbf{X}_0^d \equiv\hat{\mathbf{X}}_0^d$ and for all $s \geq1$,
$\mathbf{X}_s^d = \hat{\mathbf{X}}_{U_s^d}^d$.\vspace*{-2pt} Obviously the above
process can be reversed by setting $\hat{\mathbf{X}}_t^d$ equal to
the $t$th accepted move in the RWM process.

For each $d \geq1$, the components of $\mathbf{X}_0^d$ are
independent and identically distributed. Therefore we focus
attention on the first component as this is indicative of the
behavior of the whole process. For $d \geq1$ and $t \geq0$, let
$V_t^d = X_{[d^2 t],1}^d$ and $\hat{V}_t^d = \hat{X}_{[d^2 t],1}^d$.
%
%
\begin{theorem}
\label{main} Fix $l >0$. For all $d \geq1$, let $\mathbf{X}_0^d
\sim\pi_d$. Then, as $\dr$,
\[
V^d \Rightarrow V
\]
in the Skorokhod topology on $D[0,\infty)$, where $V_\cdot$ satisfies
the (reflected) Langevin SDE on $[0,1]$
%
%
\begin{equation}
\label{eq2a1} d V_t = \sqrt{\phi(l)} \,dB_t + \tfrac{1}{2} \phi(l)
g'(V_t) \,dt + dL^0_t(V) - dL^1_t(V)
\end{equation}
with $V_0 \sim f$. Note that $B_t$ is standard Brownian motion,
\[
\phi(l) = \frac{l^2}{3} \exp\biggl( - \frac{f^\ast l}{2} \biggr)
\]
and
$f^\ast= \lim_{x \downarrow0} ( \frac{f(x) + f(1-x)}{2}
)$.

Here $\{L^y_t, t \ge0 \}$ denotes the local time of $V$ at
$y$ $(\mbox{$=$}0,1)$ and the SDE~(\ref{eq2a1}) corresponds to standard
reflection at the boundaries $0$ and $1$ (see, e.g., Chapter
\textup{VI}
of~\cite{RevYor}).
\end{theorem}
\begin{pf}
As noted in Section~\ref{SecInt}, the acceptance
probability of the RWM process is inherently random and therefore it
is necessary to consider the behavior of the RWM process averaged
over $[d^\delta]$ iterations, for $\delta>0$. Fix $0 < 20 \gamma<
\beta< \delta< \delta+ \gamma< \frac{1}{2}$ and let $\{ k_d \}$
be a sequence of positive integers satisfying $[d^\beta] \leq k_d
\leq[d^\delta]$. For $s \in\mathbb{Z}^+$, let
$\tilde{\mathbf{X}}_s^d = \mathbf{X}_{s[d^\delta]}^d$ and for $t
\geq0$, let $\tilde{V}_t^d = \tilde{X}_{[t d^2/[d^\delta]],1}^d$.
For all $t \geq0$, $|X_{t+1,1}^d - X_{t,1}^d| \leq\sigma_d$ and
$|[d^2 t] - [d^\delta] \times[d^2 t/[d^\delta]]| \leq[d^\delta]$.
Hence, for all $T > 0$,
%
%
\begin{equation} \label{eq45c} \sup_{0 \leq s \leq T} |
\tilde{V}_s^d - V_s^d| \leq[d^\delta] \sigma_d.
\end{equation}
Therefore by~\cite{Bill}, Theorem 4.1, $V^d \Rightarrow V$ as $\dr$,
if $\tilde{V}^d \Rightarrow V$ as $\dr$. Hence we proceed by showing
that
%
%
\begin{equation}
\label{eqn112} \tilde{V}^d \Rightarrow V \qquad\mbox{as }
\dr.
\end{equation}

Let $G_d^\delta$ be the (discrete-time) generator of
$\tilde{\mathbf{X}}^d$ and let $H$ be an arbitrary test function of
the first component only. Thus
%
%
\begin{equation}
\label{eqn113} G_d^\delta H (\mathbf{x}^d) =
\frac{d^2}{[d^\delta]} \ez[ H (\tilde{\mathbf{X}}_1^d) - H
(\tilde{\mathbf{X}}_0^d) | \tilde{\mathbf{X}}_0^d = \mathbf{x}^d].
\end{equation}
The generator $G$ of the (limiting) one-dimensional diffusion $V$
for an arbitrary test function $H$ is given by
%
%
\begin{equation} \label{eqmain1}
G H(x) =\phi(l) \bigl\{ \tfrac{1}{2} g'(x) H'(x) + \tfrac{1}{2}
H''(x) \bigr\}
\end{equation}
for all $x \in[0,1]$ at least for all $H \in{\cal D}$, where
${\cal D}$ is defined in (\ref{eqmain1a}) below.

First note that the diffusion defined by (\ref{eqmain1}) is
regular; see~\cite{EK}, page 366. Therefore by~\cite{EK}, Chapter 8,
Corollary 1.2, it is sufficient to restrict attention to functions
%
%
\begin{equation} \label{eqmain1a} H \in\mathcal{D} \equiv\{ h\dvtx h
\in\hat{C} ([0,1]) \cap C^2 ((0,1)) \cap
\mathcal{D}^\ast, Gh \in\hat{C} ([0,1]) \},
\end{equation}
where
$C^2 ((0,1))$ is the set of twice differentiable functions upon
$(0,1)$, $\hat{C} [0,1]$ is the set of bounded continuous functions
upon $[0,1]$ and $\mathcal{D}^\ast$ is obtained by setting $q_i=0$
$(i=0,1)$ in~\cite{EK}, page 367, (1.11) and is given by
%
%
\begin{equation} \label{eqmain2} \mathcal{D}^\ast= \{ h\dvtx
h'(0)=h'(1) =0
\}.
\end{equation}
%
Let $H^\ast_1 = \sup_{0 \leq y \leq1} H^\prime(y)$ and $H^\ast_2 =
\sup_{0 \leq y \leq1} H^{\prime\prime} (y)$. Then $H \in C^2
((0,1))$ combined with $H \in\mathcal{D}^\ast$ implies that
$H^\ast_1 < \infty$. It then follows from $g^\prime$ being bounded
on $[0,1]$ and $GH \in\hat{C} ([0,1])$ that $H^\ast_2 < \infty$.
These observations will play a key role in Appendix~\ref{SecGen}.

Now (\ref{eqn112}) is proved using~\cite{EK}, Chapter 4, Corollary
8.7, by showing that there exists a sequence of sets $\{ \tilde{F}_d
\}$ such that for any $T >0$,
%
%
\begin{equation}
\label{eqn114} \pz\Biggl( \bigcup_{j=0}^{[T d^2/[d^\delta]]} \{
\mathbf{X}_j^d \notin\tilde{F}_d \} \Biggr) \rightarrow0
\qquad\mbox{as } \dr
\end{equation}
and
%
%
\begin{equation}
\label{eqn115}\sup_{\mathbf{x}^d \in\tilde{F}_d} | G_d^\delta H
(\mathbf{x}^d) - G H (x_1) | \rightarrow0 \qquad\mbox{as }
\dr.
\end{equation}

Let the sets $\{F_d \}$ and $\{\tilde{F}_d \}$ be such that $F_d =
\bigcap_{j=1}^4 F_d^j$ and
%
%
\begin{equation} \label{emainx1}
\tilde{F}_d = \Biggl\{ \mathbf{x}^d ; \pz\Biggl(
\bigcup_{j=0}^{[d^\delta]} \{ \hat{\mathbf{X}}_j^d \notin F_d \} |
\hat{\mathbf{X}}_0^d = \mathbf{x}^d \Biggr) \leq d^{-3} \Biggr\},
\end{equation}
where $F_d^1$, $F_d^2$, $F_d^3$ and $F_d^4$ are defined below. Recall
that $b_d^r (\mathbf{x}^d) = \sum_{j=1}^d 1_{ \{ x_j \in R_d^r
\}}$, the total number of components of $\mathbf{x}^d$ in $R_d^r
=(0,r/d) \cup(1-r/d,1)$. We term~$R_d^l$ the rejection
region, in that, for any component in $R_d^l$, there is positive
probability of proposing a move outside the hypercube with such
moves automatically being rejected. Let
%
%
\begin{eqnarray}
\label{eqn11f1}
F_d^1 &=& \{ \mathbf{x}^d ; b_d^l
(\mathbf{x}^d) \leq\gamma\log d \}, \\
\label{eqn11f2} F_d^2 &=& \bigcap_{k=[d^\beta]}^{[d^\delta]} \bigl\{
\mathbf{x}^d ; |b_d^{k^{3/4}} (\mathbf{x}^d) -
\ez[ b_d^{k^{3/4}} (\mathbf{X}_0^d)] | \leq\sqrt{k} \bigr\}, \\
\label{eqn11f3} F_d^3 &=& \Bigl\{ \mathbf{x}^d ;
\sup_{[d^\beta]\leq k_d \leq[d^\delta]} \sup_{0 \leq r \leq l} |
\lambda_d (\mathbf{x}^d;r;k_d) - \lambda(r)| \leq d^{-\gamma}
\Bigr\}, \\
\label{eqn11f4} F_d^4 &=& \Biggl\{ \mathbf{x}^d; \Biggl| \frac{1}{d}
\sum_{j=1}^d g^\prime(x_j)^2 - \ez_f [ g^\prime(X_1)^2] \Biggr| <
d^{-{1}/{8}} \Biggr\},
\end{eqnarray}
where $\lambda_d (\mathbf{x}^d; r; k_d) = \ez[
b_d^r (\mathbf{X}_{k_d}^d) | \mathbf{X}_0^d = \mathbf{x}^d]$ and
$\lambda(r) = f^\ast r (1+r/2l)$. In Appendix~\ref{SecSets}, we
prove (\ref{eqn114}) for the sets $\{\tilde{F}_d \}$ given in
(\ref{emainx1}). Note that (\ref{eqn114}) follows immediately
from Theorem~\ref{lem321}, (\ref{eqss46}) since $\mathbf{X}_0^d
\sim\pi_d$. An outline of the roles played by each $F_d^j$
$(j=1,2,3,4)$ is given below. For $\mathbf{x}^d \in F_d^1$
($\mathbf{x}^d \in F_d^2$) the total number of components in
(\textit{close to}) the rejection region are controlled. For $\mathbf{x}^d
\in F_d^3$ after $k_d$ iterations the total number and position of
the points $\{ \hat{\mathbf{X}}_{k_d}^d | \hat{\mathbf{X}}_0^d =
\mathbf{x}^d \}$ in $R_d^l$ are approximately from the stationary
distribution of $\hat{\mathbf{X}}^d$. Finally, for $\mathbf{x}^d \in
F_d^4$, $\frac{1}{d} \sum_{j=1}^d g^\prime(x_j)^2 \approx\ez_f
[g^\prime(X)^2]$; this is the key requirement for the sets\vspace*{1pt}
$\{F_d\}$ given in~\cite{RGG}, cf.~\cite{RGG}, page 114, $R_n (x_2,
\ldots, x_n)$.

The proof of (\ref{eqn115}) splits into two parts and exploits the
pseudo-RWM process. Let
%
%
\begin{equation} \label{eqmain7}
P_d = \max\Biggl\{ K=0,1, \ldots, [d^\delta-1]; \frac{1}{[d^\delta]}
\sum_{j=0}^{K-1} M_j (J_d (\hat{\mathbf{X}}_j^d)) \leq1 \Biggr\}
\Big/[d^\delta],\hspace*{-25pt}
\end{equation}
the proportion of accepted moves in the first $[d^\delta]$
iterations, where the sum is set equal to zero if vacuous. Then
$\tilde{\mathbf{X}}_1^d = \mathbf{X}_{[d^\delta]}^d =
\hat{\mathbf{X}}_{[P_d d^\delta]}^d$ and
%
%
\begin{equation}
\label{eqn116} G_d^\delta H (\mathbf{x}^d) =
\frac{d^2}{[d^\delta]} \ez\bigl[ H \bigl(\hat{\mathbf{X}}_{[P_d d^\delta]}^d\bigr)
- H (\hat{\mathbf{X}}_0^d) | \hat{\mathbf{X}}_0^d = \mathbf{x}^d\bigr].
\end{equation}
In Appendix~\ref{SecPd}, we show that for all $\mathbf{x}^d \in
\tilde{F}_d$, $P_d | \hat{\mathbf{X}}^d_0 = \mathbf{x}^d \convp
\exp
(-lf^\ast/2)$ as $\dr$. Consequently, it is useful to introduce
$\hat{G}_d^{\delta, \pi} H (\mathbf{x}^d)$ $(0 \leq\pi\leq1)$
which is defined for fixed $0 \leq\pi\leq1$ as
%
%
\begin{eqnarray} \label{eqmain6}
\hat{G}_d^{\delta, \pi} H (\mathbf{x}^d) &=& \frac{d^2}{[d^\delta]}
\ez\bigl[ \bigl(H\bigl(\hat{\mathbf{X}}^d_{[\pi d^\delta]}\bigr) - H
(\hat{\mathbf{X}}^d_0)\bigr) | \hat{\mathbf{X}}^d_0 = \mathbf{x}^d
\bigr] \nonumber\\
&=& \frac{d^2}{[d^\delta]} \sum_{j=0}^{[\pi d^\delta-1]} \ez
[ H (\xx_{j+1}^d) - H (\xx_j^d) | \xx_0^d = \mathbf{x}^d
]
\\
&=& \frac{1}{[d^\delta]} \sum_{j=0}^{[\pi d^\delta-1]} \ez[
\hat{G}_d H (\xx_j^d) | \xx_0^d = \mathbf{x}^d ],\nonumber
\end{eqnarray}
where
%
%
\begin{equation} \label{eqmain6x} \hat{G}_d H (\xx_j^d) = d^2 \ez
[ H (\xx_1^d - \xx_0^d) | \xx_0^d
= \mathbf{x}^d].
\end{equation}

Finally in Appendix~\ref{SecGen}, we prove in Lemma~\ref{lem45} that
%
%
\begin{equation} \label{eqn117}
\sup_{0 \leq\pi\leq1} \sup_{\mathbf{x}^d \in\tilde{F}_d}
|\hat{G}_d^{\delta, \pi} H (\mathbf{x}^d) - G H (x_1)| \rightarrow0
\qquad\mbox{as } \dr.
\end{equation}
The triangle
inequality is then utilized to prove (\ref{eqn115}) in Lemma~\ref{lem45}
using~(\ref{eqn117}) and $P_d | \mathbf{X}^d_0 =
\mathbf{x}^d \convp\exp(-lf^\ast/2)$ as $\dr$.
\end{pf}

It should be noted that in Appendix~\ref{SecGen}, we assume that $\ez
[g^\prime(X)^2] > 0$, in particular in Lemma~\ref{lem40}. In
Appendixes~\ref{SecSets} and~\ref{SecPd} we make no such assumption.
However, $\ez[g^\prime(X)^2]=0$ corresponds to $f (x) =1_{\{ 0 < x
< 1 \}}$ (uniform distribution), and proving Lemma~\ref{lem45} in
this case follows similar but simpler arguments to those given in
Appendix~\ref{SecGen}.

A key difference between the diffusion limits for continuous and
discontinuous i.i.d. product densities is the dependence of the
speed measure $\phi(l)$ upon $f$. For continuous (suitably
differentiable) $f$, $\phi(l)$ depends upon $I \equiv\ez_f [ \{
f^\prime(X)/f (X) \}^2]$, which is a measure of the ``roughness'' of~$f$.
For discontinuous densities of the form (\ref{eq1b}), $\phi
(l)$ depends upon $f^\ast= \lim_{x \downarrow0} \{ (f(x) +
f(1-x))/2 \}$, the (mean of the) limit of the density at the
boundaries (discontinuities). Discussion of the role of the density
$f$ in the behavior of the RWM algorithm is given in Section
\ref{SecExtSim}.

The most important consequence of Theorem~\ref{main} is the
following result.
%
%
\begin{coro}
\label{mc1} Let $a(l) = \exp(-f^\ast l/2)$. Then
\[
\ez_{\pi_d} \ez[ J_d (\mathbf{X}_0^d) ] \rightarrow
a(l) \qquad\mbox{as } \dr.
\]
$\phi(l)$ is maximized by $l= \hat{l} = 4/ f^\ast$ with
\[
a(\hat{l}) = \exp(-2) = 0.1353.
\]
\end{coro}

Clearly, if $f (\cdot)$ is known, $\hat{l}$ can be calculated
explicitly. However, where MCMC is used, $f (\cdot)$ will often only
be known up to the constant of proportionality. This is where
Corollary~\ref{mc1} has major practical implications, in that, to
maximize the speed of the limiting diffusion, and hence, the
efficiency of the RWM algorithm, it is sufficient to monitor the
average acceptance rate, and to choose $l$ such that the average
acceptance rate is approximately $e^{-2}$. Therefore there is no
need to explicitly calculate or estimate the constant of
proportionality.

\section{Extensions} \label{SecExtSim}

In this section, we discuss the extent to which the conclusions of
Theorem~\ref{main} extend beyond $\pi_d$ being an i.i.d. product
density upon the $d$-dimensional hypercube and $\mathbf{X}_0^d \sim
\pi_d$. First we present two extensions of Theorem~\ref{main}. The
second extension, Theorem~\ref{thmprop}, is an important practical
result concerning lower-dimensional updating schema.

Suppose that $f (\cdot)$ is nonzero on the positive half-line. That
is,
%
%
\begin{equation} \label{eqtri1} f (x) \propto\exp(g(x))
\qquad(x > 0)
\end{equation}
and $f(x) =0$ otherwise.
%
%
\begin{theorem}
\label{thmhalf} Fix $l >0$. For all $d \geq1$, let $\mathbf{X}_0^d
\sim\pi_d$, given by (\ref{eqtri1}), with $\sup_{x \geq0}
|g^\prime(x)| = g^\ast< \infty$. Then, as $\dr$,
\[
V^d \Rightarrow V
\]
in the Skorokhod topology on $D[0,\infty)$, where $V_\cdot$ satisfies
the (reflected) Lan\-gevin SDE on $[0, \infty)$
\[
d V_t =
\sqrt{\phi(l)} \,dB_t + \tfrac{1}{2} \phi(l) g'(V_t) \,dt + dL_t^0
(V)
\]
with $V_0 \sim f$,
$\phi(l) = \frac{l^2}{3} \exp(-f^\star l/4)$ and $f^\star=
\lim_{x \downarrow0} f (x)$.
\end{theorem}
\begin{pf}
The proof of the theorem is virtually identical to the
proof of Theorem~\ref{main}, and so, the details are omitted.
\end{pf}

Note that we have assumed that $g^\prime(\cdot)$ is bounded on $[0,
\infty)$. This assumption is almost certainly stronger than
necessary with $g^\prime(\cdot)$ being Lipschitz and/or satisfying
certain moment conditions probably being sufficient; cf.~\cite{RGG}.

Theorem~\ref{thmhalf} is unsurprising with the speed of the
diffusion depending upon the number of components close to the
discontinuity at 0.
%
%
\begin{coro}
\label{mc2} Let $\pi_d (\mathbf{x}^d) = \prod_{i=1}^d f (x_i)$ where
$f$ satisfies (\ref{eqtri1}). Then
\[
\ez_{\pi_d} [ J_d (\mathbf{X}_0^d) ] \rightarrow\exp
(-f^\star l/4) \equiv a(l) \qquad\mbox{as } \dr.
\]
$\phi(l)$ is maximized by $l= \hat{l} = 8/ f^\star$ with
\[
a(\hat{l}) = \exp(-2) = 0.1353.
\]
\end{coro}

Therefore the conclusions are identical to Corollary~\ref{mc1} that
in order to maximize the speed of the limiting diffusion it is
sufficient to choose $l$ such that the average acceptance rate is
$e^{-2}$.

The second and more important extension of Theorem~\ref{main}
follows on from~\cite{NR06}. In~\cite{NR06}, the
Metropolis-within-Gibbs algorithm was considered, where only a
proportion $c$ $(0 < c \leq1)$ of the components are updated at
each iteration. For given $d \geq1$, at each iteration $c_d d$ of
the components are chosen uniformly at random and new values for
these components are proposed using random walk Metropolis with
proposal variance $\sigma_{d,c_d}^2 =(l/d)^2$. The remaining
$(1-c_d) d$ components remain fixed at their current values.
Finally, it is assumed that $c_d \rightarrow c$ as $\dr$.

The following result assumes that $f (\cdot)$ is nonzero on $(0,1)$
only. The extension to the positive half-line is trivial.
%
%
\begin{theorem}
\label{thmprop} Fix $0 < c \leq1$ and $l >0$. For all $d \geq1$,
let $\mathbf{X}_0^d = (X_{0,1}^d, X_{0,2}^d,\break \ldots, X_{0,d}^d)$ be
such that all of its components are distributed according to $f
(\cdot)$. Then, as $\dr$,
\[
V^d \Rightarrow V
\]
in the Skorokhod topology, where $V_0 \sim f (\cdot)$ and $V$ satisfies
the (reflected) Langevin SDE on $[0, 1]$
\[
d V_t = \sqrt{\phi_c (l)} \,dB_t + \tfrac{1}{2} \phi_c
(l) g'(V_t) \,dt + dL_t^0 (V) - dL_t^1 (V),
\]
where $B_t$ is standard Brownian motion,
$\phi_c (l) = \frac{c l^2}{3} \exp(-c f^\ast l/2 )$ and $f^\ast=
\lim_{x \downarrow0} \frac{f (x) + f (1-x)}{2}$.
\end{theorem}

Let $a_d^{c_d} (l)$ denote the average acceptance rate of the RWM
algorithm in~$d$ dimensions where a proportion $c_d$ of the
components are updated at each iteration. Let
\[
a^c(l) = \exp(- cf^\ast l/2).
\]
We then have the following
result which mirrors Corollaries~\ref{mc1} and~\ref{mc2}.
%
%
\begin{coro}
\label{mc3} Let $c_d \rightarrow c$ as $\dr$. Then
\[
a_d^{c_d} (l) \rightarrow a^c(l) \qquad\mbox{as } \dr.
\]
For fixed $0 < c \leq1$, $\phi_c (l)$ is maximized by
\[
l= \hat{l}_c = \frac{4}{c f^\ast}
\]
and
\[
\phi_c (\hat{l}_c) = \frac{1}{c} \phi_1 (\hat{l}_1).
\]
Also
\[
a(\hat{l}_c) = \exp(-2) = 0.1353.
\]
\end{coro}

Corollary~\ref{mc3} is of fundamental importance from a practical
point of view, in that it shows that the optimal speed of the
limiting diffusion is inversely proportional to $c$. Therefore the
optimal action is to choose $c$ as close to 0 as possible.
Furthermore, we have shown that not only is full-dimensional RWM bad
for discontinuous target densities but it is the worst algorithm of
all the Metropolis-within-Gibbs RWM algorithms.

We now go beyond i.i.d. product densities with a discontinuity at
the boundary and $\mathbf{X}_0^d \sim\pi_d$. We consider general
densities on the unit hypercube, discontinuities not at the boundary
and $\mathbf{X}_0^d \not\sim\pi_d$. As mentioned in Section
\ref{SecInt}, for i.i.d. product densities, the speed measure of
the limiting one-dimensional diffusion, $\phi(l)$, is equal to the
limit, as $\dr$, of the ESJD times $d$. Therefore we consider the
ESJD for the above-mentioned extensions as being indicative of the
behavior of the limiting Langevin diffusion. We also highlight an
extra criterion which is likely to be required in moving from an
ESJD to a Langevin diffusion limit.

Using the proof of Theorem~\ref{main}, it is straightforward to show
that
%
%
\begin{eqnarray} \label{eqextd1}
\phi(l) &=& \frac{l^2}{3} \exp\biggl( - \frac{l f^\ast}{2}
\biggr) \nonumber\\
&=& l^2 \ez[ Z_1^2] \lim_{\dr} \ez\bigl[ 1_{ \{ \mathbf{X}_0^d +
\sigma_d \mathbf{Z}_1^d \in[0,1]^d \}} \bigr]
\\
&=& \frac{l^2}{3} \lim_{\dr} \ez\biggl[ \biggl( \frac{3}{4}
\biggr)^{b_d^l (\mathbf{X}_0^d)} \biggr].\nonumber
\end{eqnarray}
The first
equality in (\ref{eqextd1}) can be proved using Lemma~\ref{lem33},
(\ref{eq33b}), where for $Z_1 \sim U(-1,1)$, $\ez[Z_1^2] = 1/3$. The
second equality in (\ref{eqextd1}) comes from
the fact that for $0 < x < \sigma_d$, $f (x) + f(1-x) = 2 f^\ast+
O(1/d)$ and for a component $X_{0,i}^d$ uniformly distributed on
$(0, l/d)$ or $(1-l/d,1)$, $P(X_{0,i}^d + \sigma_d Z_i^d \in[0,1])
=3/4$. That is, the acceptance probability of a proposed move is
dominated by whether or not the proposed move lies inside the
$d$-dimensional unit hypercube. Proposed moves inside the hypercube
are accepted with probability $1 - o(d^{-\alpha})$ for any $\alpha<
1/2$; see Lemma~\ref{lem35}. Thus it is the number and behavior of
the components at the boundary of the hypercube (the discontinuity)
which determine the behavior of the RWM algorithm. This is also
seen in Theorems~\ref{thmhalf} and~\ref{thmprop}.

First, we consider discontinuities not at the boundary. Suppose
that\break $\pi_d (\mathbf{x}^d) = \prod_{i=1}^d f (x_i)$, where
%
%
\begin{equation} \label{eqextb1}
f(x) \propto1_{ \{ x \in[a,b] \}} \exp(g (x))
\qquad(x \in\mathbb{R})
\end{equation}
for some $a, b \in\mathbb{R}$.
Further suppose that $g (\cdot)$ is continuous (twice
differentiable) upon $[a, b]$ except at a countable number of points,
$\mathcal{P} = \{ a_1, a_2, \ldots, a_k \}$, say, on $(a,b)$. Set
$a_0 =a$ and $a_{k+1} =b$, with $\sigma_d = l/d$. For $j=0,1,
\ldots, k+1$, let $f_j^- = \lim_{x \rightarrow a_j -} f (x)$ and
$f_j^+ = \lim_{x \rightarrow a_j +} f (x)$, with $Y_j^- \sim
\operatorname{Po} (l f_j^-/4)$ and $Y_j^+ \sim\operatorname{Po} (l
f_j^+/4)$, where
$f_0^-=f_{k+1}^+=Y_0^-=Y_{k+1}^+=0$. Then\vspace*{1pt} following~\cite{NR10},
(4.23), we can show that $d$ times the ESJD
%
%
\begin{equation}
\label{eqem1} d \ez\Biggl[ \sum_{i=1}^d (X_{1,i}^d- X_{0,i}^d)^2
\Biggr] \rightarrow\frac{l^2}{3} \ez\Biggl[ 1 \wedge
\prod_{j=0}^{k+1} \biggl( \frac{f_j^-}{f_j^+} \biggr)^{Y_j^+ - Y_j^-}
\Biggr] \qquad\mbox{as } \dr.\hspace*{-25pt}
\end{equation}
Thus the
optimal scaling of $\sigma_d$ is again of the form $l/d$ and the
acceptance or rejection of a proposed move is determined by the
components close to the discontinuities. Furthermore, it is
straightforward to show that for each $j=0, 1, \ldots, k+1$, $l^2
(f_j^-/f_j^+)^{Y_j^+ - Y_j^-} \convp0$ as $l \rightarrow\infty$,
implying\vspace*{1pt} that the optimal choice of $l$ lies in $(0, \infty)$.
Proving a Langevin diffusion for the (normalized) first component of
the RWM algorithm should be possible with appropriate local time
terms at the discontinuities in $f$. While (\ref{eqem1}) holds
regardless of $f_j^-$ and $f_j^+$ for a diffusion limit we require
that\vspace*{1pt} $\min_{1 \leq j \leq k+1} f_j^-, \min_{0 \leq j \leq k} f_j^+
>0$, that is, the density is strictly positive on $(a,b)$. (If this
is not the case, the RWM algorithm is reducible in the limit as
$\dr$.) Extensions to the case where either $a=-\infty$ and/or $b=
\infty$ are straightforward.

Second, we consider general densities which are zero outside the
$d$-dimen\-sional hypercube, $\pi_d (\mathbf{x}^d) \propto1_{ \{
\mathbf{x}^d \in[0,1]^d \}} \exp(\mu_d (\mathbf{x}^d))$, where
$\mu_d (\cdot)$ is assumed to be continuous and twice
differentiable. Let $\sigma_d =l/d$ and assuming that
%
%
\begin{equation} \label{eqextd2} \exp\bigl( \mu_d (\mathbf{X}_0^d
+ \sigma_d \mathbf{Z}_1^d) - \mu_d (\mathbf{X}_0^d) \bigr) \convp1
\qquad\mbox{as } \dr,
\end{equation}
we have that $d$ times the ESJD satisfies
%
%
\begin{equation} \label{eqextd3} d \ez\Biggl[ \sum_{i=1}^d
(X_{1,i}^d - X_{0,i}^d)^2 \Biggr] \rightarrow\frac{l^2}{3}
\lim_{\dr} \ez\biggl[ \biggl( \frac{3}{4} \biggr)^{b_d^l
(\mathbf{X}_0^d)} \biggr] \qquad\mbox{as } \dr.
\end{equation}
Note that (\ref{eqextd2}) is a weak condition and should be
straightforward to check using a Taylor series expansion of $\mu_d$.
For\vspace*{1pt} i.i.d. product densities, $b_d^l (\mathbf{X}_0^d ) \convd
\operatorname{Po} (2 l f^\ast)$ as $\dr$. More generally, the limiting
distribution of $b_d^l (\mathbf{X}_0^d)$ will determine the limit of
the right-hand side of (\ref{eqextd3}). In particular, so long as
there exist $\delta>0$ and $K \in\mathbb{N}$ such that $\pz
(\lim_{\dr} b_d^l (\mathbf{X}_0^d) \leq K) > \delta$, the right-hand
side of (\ref{eqextd3}) will be nonzero for $l >0$. It is
informative to consider what conditions upon $\pi_d$ are likely to
be necessary for a diffusion limit, whether it be one-dimensional or
infinite-dimensional as in~\cite{BR00}. Suppose that \mbox{$b_d^l
(\mathbf{X}_0^d) \convd B$} as $\dr$. For a diffusion limit we will
require moment conditions on $B$, probably requiring that there
exists $\varepsilon
>0$ such that $\ez[ \exp(\varepsilon B)] < \infty$. This will be
required to control the probability of the RWM algorithm getting
``stuck'' in the corners of the hypercube. This highlights a key
difference between studying the ESJD and a diffusion limit. For the
ESJD,\vadjust{\goodbreak} we want a~positive probability that the total number of
components at the boundary of the hypercube is finite in the limit
as $\dr$. For the diffusion limit, as seen with the construction of
$\{F_d^1 \}$ in Theorem~\ref{main}, we want that the probability of
there being a large number of components $(O(\log d))$ at the
boundary is very small $(o(d^{-2}))$.

Third, suppose that $\mathbf{X}_0^d \not\sim\pi_d$. There are
very bad starting points in the ``corners'' of the hypercube. For
example, if $\mathbf{X}_0^d = (\exp(-d), \exp(-d), \ldots,\break
\exp(-d))$, $J_d (\mathbf{X}_0^d) \approx(0.5 +\exp(-d))^d$ which
even for $d=100$ is less than $1 \times10^{-30}$. Thus the RWM
process is likely to be ``stuck'' at its starting point for a very
long period of time. This is rather pathological and a more
interesting question is the situation when $\mathbf{X}_0^d
=\mathbf{S}^d$, where the components of $\mathbf{S}^d$ are i.i.d. In
particular, suppose that $S_1^d \sim U[0,1]$, so that
$\mathbf{X}_0^d$ is chosen uniformly at random over the hypercube.
Note that, if $\mathbf{S}^d$ is the uniform distribution,
%
%
\begin{equation}
\label{eqrev1} d \ez\Biggl[ \sum_{i=1}^d (X_{1,i}^d - X_{0,i}^d)^2
| \mathbf{X}_0^d \eqd\mathbf{S}^d \Biggr] \rightarrow\frac{l^2}{3}
\exp\biggl( - \frac{l}{2} \biggr) \qquad\mbox{as } \dr
\end{equation}
with the right-hand side maximized by taking $\hat{l}
=4$ compared with $\hat{l} = 4/f^\ast$ for $\mathbf{X}_0^d \sim
\pi_d$. We expect to see similar behavior to~\cite{CRR05}, in that
the optimal~$\sigma_d$ (in terms of the ESJD) will vary as the
algorithm converges to the stationary distribution but will be of
the form $\sigma_d = l/d$ throughout. The RWM algorithm is unlikely
to get ``stuck'' with it conjectured that for any $T >0$ and $\gamma
>0$,
\[
\pz\Biggl( \bigcup_{t=0}^{[T d^2]} \{ b_d^l (\mathbf{X}_t^d) \geq
\gamma\log d \}
| \mathbf{X}_0^d \eqd\mathbf{S}^d \Biggr) \rightarrow0
\qquad\mbox{as } \dr.
\]
Simulations with
$f (x) \propto1_{ \{ 0 < x < 1 \}} \exp(-2x)$ and $f(x) \propto
1_{ \{ 0 < x < 1 \}} \exp(-(x-0.5)^2/2)$ and $d=10, 20, \ldots,
200$ suggest that convergence occurs in $O(d^2)$ \mbox{iteration}. For
convergence, we monitor the mean of $\mathbf{X}_t^d$ for $f (x)
\propto\break 1_{ \{ 0 < x < 1 \}} \exp(-2x)$ and the variance of
$\mathbf{X}_t^d$ for $f (x) \propto1_{ \{ 0 < x < 1 \}} \exp
(-(x-0.5)^2/2)$.

\begin{appendix}\label{app}

\section{\texorpdfstring{Construction of the sets $\{F_d\}$ and $\{\tilde{F}_d\}$}
{Construction of the sets $\{F_d\}$ and $\{F_d\}$}}
\label{SecSets}

The sets $F_d$ consist of the intersection of four sets $F_d^i$
$(i=1,2,3,4)$. For $i=1,2,3,4$, we will define $F_d^i$ and discuss
the role that it plays in the proof of Theorem~\ref{main}, one at a
time. Furthermore, we show that in stationarity it is highly
unlikely that $\mathbf{X}_t^d$ does not belong to $F_d$. Since we
rely upon a~homogenization argument, it is necessary to go further
than the sets $F_d$ to the sets $\tilde{F}_d \subset F_d$. In
particular, if $\hat{\mathbf{X}}_0^d \in\tilde{F}_d$, then it is
highly unlikely that any of $\hat{\mathbf{X}}_1^d,
\hat{\mathbf{X}}_2^d, \ldots,\hat{\mathbf{X}}_{[d^\delta]}^d$ do not
belong to $F_d$. The above statement is made precise in Theorem
\ref{lem321} below, where the constructions of $\{F_d \}$ and $\{
\tilde{F}_d \}$ are drawn together.\vadjust{\goodbreak}

It is possible that all $d$ components of $\mathbf{X}_0^d$ are in
$R_d^l$. However, this is highly unlikely and we show in Lemma
\ref{lem31} that with high probability, there are at most
$\gamma\log d$ components in the rejection region. Let $F_d^1 = \{
\mathbf{x}^d ; b_d^l (\mathbf{x}^d) \leq\gamma\log d \}$.\vspace*{-1pt}
%
%
\begin{lem}
\label{lem31} For any $\kappa> 0$,
\[
d^\kappa\pz(\mathbf{X}_0^d
\notin F_d^1) \rightarrow0 \qquad\mbox{as } \dr.\vspace*{-1pt}
\]
\end{lem}
\begin{pf}
Fix $ \kappa> 0$. Note that $\mathbf{X}_0^d \notin
F_d^1$ if and only if $b_d^l (\mathbf{X}_0^d) > \gamma\log
d$. However,
\[
b_d^l (\mathbf{X}_0^d) \sim\operatorname{Bin} \biggl(d, \int_0^{l/d} \{ f(x)
+ f
(1-x) \} \,dx \biggr)\vspace*{-2pt}
\]
with
%
%
\begin{equation} \label{eq31a}
d \int_0^{l/d} \{ f(x) + f (1-x)\} \,dx \rightarrow2 f^\ast l
\qquad\mbox{as } \dr.
\end{equation}

Fix $\rho> \kappa/\gamma$. By Markov's inequality and using
independence of the components of $\mathbf{X}_0^d$,
%
%
\begin{eqnarray} \label{eq31b}
&&
d^\kappa\pz\bigl(b_d^l (\mathbf{X}_0^d) > \gamma\log d\bigr)\nonumber\\[-1pt]
&&\qquad
\leq d^\kappa\ez[ \exp(\rho b_d^l (\mathbf{X}_0^d))]/ \exp
(\rho\gamma\log d) \nonumber\\[-1pt]
&&\qquad= d^\kappa\ez\bigl[ \exp\bigl(\rho1_{\{ X_{0,1}^d \in R_d^l
\}}\bigr)\bigr]^d / d^{\rho\gamma}
\\[-2pt]
&&\qquad = d^\kappa\biggl( 1 + (e^\rho-1) \int_0^{l/d} \{ f(x) +
f(1-x)\} \,dx \biggr)^d \big/ d^{\rho\gamma} \nonumber\\[-2pt]
&&\qquad \leq d^{\kappa- \rho\gamma} \exp\biggl( (e^\rho-1)d
\int_0^{l/d} \{ f(x) + f (1-x)\} \,dx \biggr).\nonumber
\end{eqnarray}
The
lemma follows since (\ref{eq31a}) implies that the right-hand side
of (\ref{eq31b}) converges to 0 as $\dr$.\vspace*{-2pt}
\end{pf}

For $\mathbf{x}^d \in F_d^1$, it follows from (\ref{eq21ay}) that
%
%
\begin{equation} \label{eqlowerJ} \quad
J_d (\mathbf{x}^d) \geq\exp
(-lg^\ast) 2^{-b_d^l (\mathbf{x}^d)}
\geq\exp(-lg^\ast) 2^{- \gamma\log d} \geq\exp(-lg^\ast)
d^{-\gamma}.
\end{equation}
This is a useful lower bound for the
acceptance probability and as a result the random walk Metropolis
algorithm does not get ``stuck'' at values of $\mathbf{x}^d \in
F_d^1$. To assist with the homogenizing arguments, we define $\{
\tilde{F}_d^1 \}$ by
%
%
\begin{equation} \label{eqtildef11}
\tilde{F}_d^1 = \Biggl\{ \mathbf{x}^d ; \pz\Biggl(
\bigcup_{j=0}^{[d^\delta]} \hat{\mathbf{X}}_j^d \notin F_d^1 |
\hat{\mathbf{X}}_0^d = \mathbf{x}^d \Biggr) \leq d^{-3} \Biggr\}.
\end{equation}
That is, by starting in $\tilde{F}_d^1$ it is highly unlikely that
the pseudo-RWM algorithm leaves $F_d^1$ in $[d^\delta]$ iterations.
To study $\tilde{F}_d^1$ and later $\tilde{F}_d$ we require the
following lemmas.\vadjust{\goodbreak}
%
%
\begin{lem} \label{lem31a} For a random variable $X$, suppose that
there exist
$\delta, \varepsilon
> 0$ such that
%
%
\begin{equation} \label{eq31c} \pz(X \in A | X \in B) \leq\delta
\varepsilon
\end{equation}
and for all $x \in D^C$, $\pz(X \in A |
X=x) \geq\varepsilon$. Then
%
%
\begin{equation} \label{eq31d} \pz(X \notin D | X \in B) \leq
\delta.
\end{equation}
\end{lem}
\begin{pf}
First note that
%
%
\begin{eqnarray} \label{eq31e}\quad
\pz(X\in A| X \in B) & \geq& \pz(X \in A \cap X \in D^C | X \in B)
\nonumber\\[-8pt]\\[-8pt]
&=& \pz(X \in A | X \in D^C, X \in B) \pz(X \notin D
| X \in B).\nonumber
\end{eqnarray}
The lemma follows from rearranging
(\ref{eq31e}) and using (\ref{eq31c}) and $\pz(X \in A | X \in D^C,
X \in B) \geq\varepsilon$.
\end{pf}
%
%
\begin{lem}
\label{lem34} Suppose that a sequence of sets $\{F_d^\star\}$ is
such that there exists $\kappa> 0$ such that
%
%
\begin{equation} \label{eq34a} d^\kappa\pz(\mathbf{X}_0^d \notin
F_d^\star) \rightarrow0 \qquad\mbox{ as } \dr.
\end{equation}

Fix $\varepsilon>0$ and let
%
%
\begin{equation} \label{eq34b} \tilde{F}_d^\star= \Biggl\{ \mathbf
{x}^d; \pz \Biggl( \bigcup_{i=0}^{[d^\delta]} \{ \hat{\mathbf{X}}_i^d
\notin F_d^\star\cap F_d^1 \} | \hat{\mathbf{X}}_0^d = \mathbf{x}^d
\Biggr) \leq d^{-\varepsilon} \Biggr\}.
\end{equation}
Then
%
%
\begin{equation} \label{eq34c} d^{\kappa- (2+ \delta+ \gamma+
\varepsilon)} \pz(\hat{\mathbf{X}}_0^d \notin
\tilde{F}_d^\star) \rightarrow0 \qquad\mbox{ as } \dr.
\end{equation}
\end{lem}
\begin{pf}
Since $\mathbf{X}_i^d \sim\pi_d$,
%
%
\begin{equation} \label{eq34e} \pz\Biggl( \bigcup_{i=0}^{[d^{2 +
\delta+ \gamma}]} \{ \mathbf{X}_i^d
\notin F_d^\star\cap F_d^1 \} \Biggr) \leq d^{2 + \delta+
\gamma} \pz(\mathbf{X}_0^d \notin F_d^\ast\cap F_d^1).
\end{equation}
Therefore for all sufficiently large $d$,
%
%
\begin{equation} \label{eq34ea} \pz(\mathbf{X}_0^d \in F_d^\ast\cap
F_d^1) \geq1 - \pz(\mathbf{X}_0^d \notin F_d^\ast) - \pz(\mathbf
{X}_0^d \notin F_d^1) \geq\tfrac{1}{2}.
\end{equation}
By Bayes's theorem, $\pz(A |B) = \pz(A \cap B)/\pz(B) \leq\pz
(A)/\pz(B)$. Therefore taking $A = \bigcup_{i=0}^{[d^{2 + \delta+
\gamma}]} \{ \mathbf{X}_i^d \notin F_d^\ast\cap F_d^1 \}$ and $B =
\{ \mathbf{X}_0^d \in F_d^\ast\cap F_d^1 \}$, it follows from
(\ref{eq34e}) and (\ref{eq34ea}) that
%
%
\begin{equation} \label{eq34f}
\pz\Biggl( \bigcup_{i=0}^{[d^{2 +
\delta+ \gamma}]} \!\{ \mathbf{X}_i^d
\notin F_d^\star\cap F_d^1\} | \mathbf{X}_0^d \in F_d^\star\cap
F_d^1 \!\Biggr) \leq\frac{d^{2 + \delta+ \gamma} \pz(\mathbf{X}_0^d
\notin F_d^\star\cap F_d^1 )}{1/2}.\hspace*{-40pt}
\end{equation}

Let
%
%
\begin{equation} \label{eq34d} \hat{F}_d^\star= \Biggl\{ \mathbf
{x}^d; \pz
\Biggl( \bigcup_{i=0}^{[d^{2 + \delta+ \gamma}]} \{ \mathbf{X}_i^d
\notin F_d^\star\cap F_d^1 \} | \mathbf{X}_0^d = \mathbf{x}^d
\Biggr) \leq d^{-\varepsilon} \Biggr\}.
\end{equation}
It follows from
Lemmas~\ref{lem31} and~\ref{lem31a} that
%
%
\begin{equation} \label{eq34g} d^{\kappa- (2 + \delta+ \gamma+
\varepsilon)} \pz( \mathbf{X}_0^d \notin
\hat{F}_d^\star| \mathbf{X}_0^d \in F_d^\star\cap F_d^1 )
\rightarrow0 \qquad\mbox{as } \dr.
\end{equation}
Since $d^\kappa\pz(\mathbf{X}_0^d \notin F_d^\star
\cap F_d^1 ) \rightarrow0$ as $\dr$, it follows from (\ref{eq34g})
that
%
%
\begin{equation} \label{eq34h} d^{\kappa- (2 + \delta+ \gamma+
\varepsilon)} \pz( \mathbf{X}_0^d \notin
\hat{F}_d^\star) \rightarrow0 \qquad\mbox{as } \dr.
\end{equation}

For $d \geq1$ and $i=0,1,2,\ldots,$ let $\{\theta_i^d \}$ be
independent and identically distributed Bernoulli random variables
with $\pz(\theta_0^d =1) = \exp(-l g^\ast) 2^{-\gamma\log d}$
where $g^\ast= {\max_{\{0 \leq x \leq1 \}}} |g^\prime(x)|$. It is
straightforward using Hoeffding's inequality to show that
%
%
\begin{equation} \label{eq34j} d^\kappa\pz\Biggl(\sum_{i=1}^{[d^{2 +
\delta+ \gamma}]} \theta_i^d < d^\delta\Biggr) \rightarrow0
\qquad\mbox{as } \dr.
\end{equation}

Now $\{\theta_j^d \}$ and $\{ \mathbf{X}_j^d \}$ can be constructed
upon a common probability space such that if $\theta_j^d =1$ and
$\mathbf{X}_j^d \in F_d^1$, $\mathbf{X}_{j+1}^d \neq
\mathbf{X}_j^d$. For $k, n \geq0$, consider $\hat{\mathbf{X}}_k^d$,
if $\sum_{i=1}^n \theta_i^d \geq k$ and $\bigcap_{j=0}^n \{ \mathbf
{X}_j^d \in F_d^\ast\cap F_d^1 \}$, a coupling exists
such that there exists $J_k \in\{k, k+1, \ldots, n \}$ such that
$\hat{\mathbf{X}}_k^d = \mathbf{X}_{J_k}^d \in F_d^\ast\cap F_d^1$.
Exploiting the above coupling, $\bigcap_{j=0}^{[d^{2 + \delta+
\gamma}]} \{ \mathbf{X}_j^d \in F_d^\star\cap F_d^1 \}$ and
$\sum_{i=1}^{[d^{2 + \delta+ \gamma}]} \theta_i^d \geq d^\delta$
together imply that $\bigcap_{j=0}^{[d^{\delta}]} \{
\hat{\mathbf{X}}_j^d \in F_d^\star\cap F_d^1\}$. Thus
%
%
\begin{equation} \label{eq34k}
\pz(\hat{\mathbf{X}}_0^d \notin\tilde{F}_d^\star) \leq\pz
(\mathbf{X}_0^d \notin\hat{F}_d^\star)
+ \pz\Biggl(\sum_{i=1}^{[d^{2 + \delta+ \gamma}]} \theta_i^d <
d^\delta\Biggr),
\end{equation}
and (\ref{eq34c}) follows from
(\ref{eq34h}), (\ref{eq34j}) and (\ref{eq34k}).
\end{pf}

As noted in Section~\ref{SecAlg}, we follow~\cite{NR08} by
considering the behavior of the random walk Metropolis algorithm
over steps of size $[d^\delta]$ iterations. We find that a single
component moves only a small distance in $[d^\delta]$ iterations,
while over $[d^\delta]$ iterations the acceptance probability,
which is dominated by the number and position of components in
$R_d^l$, ``forgets'' its starting value. Moreover, we show that
approximately $\exp(-f^\ast l/2) [d^\delta]$ of the proposed moves
are accepted. However, we need to control the number of components
which are \textit{close to} the rejection region ($F_d^2$) and the
distribution of the position of the components in the rejection
region after $[d^\beta]$ iterations ($F_d^3$), where $0 < \beta<
\delta$.

For any $k \geq1$, let
\[
\hat{F}_d^2 (k) = \bigl\{ \mathbf{x}^d\dvtx
|b_d^{k^{3/4}} (\mathbf{x}^d) - \ez[ b_d^{k^{3/4}}
(\mathbf{X}^d_0)] | \leq\sqrt{k} \bigr\}
\]
and let
%
%
\begin{equation} \label{eqssz1}
F_d^2 = \bigcap_{k=[d^\beta]}^{[d^\delta]} \hat{F}_d^2 (k).
\end{equation}

Before studying $F_d^2$, we state a simple, useful result concerning
the central moments of a sequence of binomial random variables.
%
%
\begin{lem} \label{lemz32} Let $B_d \sim\operatorname{Bin} (d, p_d)$.
Suppose that
$p_d \rightarrow0$ and $d p_d
\rightarrow\infty$ as $\dr$; then for any $m \in\mathbb{N}$,
%
%
\begin{equation} \label{eqlz321} \ez\bigl[(B_d - \ez
[B_d])^{2m}\bigr]/(dp_d)^m \rightarrow\prod_{j=1}^m (2j-1)
\qquad\mbox{as } \dr.
\end{equation}
\end{lem}
%
%
\begin{lem} \label{lem32} For any $\kappa> 0$ and sequence of positive
integers $\{ k_d \}$ satisfying $[d^\beta] \leq k_d \leq
[d^\delta]$,
%
%
\begin{equation} \label{eq32c} d^\kappa\pz
\bigl(\mathbf{X}_0^d \notin\hat{F}_d^2 (k_d) \bigr) \rightarrow0
\qquad\mbox{as } \dr.
\end{equation}

Consequently, for any $\kappa> 0$, $d^\kappa\pz(\mathbf{X}_0^d
\notin F_d^2 ) \rightarrow0$ as $\dr$.
\end{lem}
\begin{pf}
Fix $\kappa> 0$. By stationarity and Markov's
inequality, for all \mbox{$m \in\mathbb{N}$},
%
%
\begin{eqnarray} \label{eq32a}
d^\kappa\pz\bigl(\mathbf{X}_0^d \notin\hat{F}_d^2 (k_d) \bigr) &=&
d^\kappa\pz\bigl( | b_d^{k_d^{3/4}} (\mathbf{X}^d_0) - \ez[
b_d^{k_d^{3/4}} (\mathbf{X}^d_0)] | \geq\sqrt{k_d}\bigr)
\nonumber\\[-8pt]\\[-8pt]
&\leq& \frac{d^\kappa}{k_d^m} \ez\bigl[ \bigl( b_d^{k_d^{3/4}}
(\mathbf{X}^d_0)
- \ez[ b_d^{k_d^{3/4}} (\mathbf{X}^d_0)]\bigr)^{2m} \bigr].\nonumber
\end{eqnarray}
However, $b_d^{k_d^{3/4}} (\mathbf{X}^d_0) \sim\operatorname{Bin} (d,
\int_0^{k_d^{3/4}/d} \{f (x) + f (1-x) \} \,dx )$, so by
Lemma~\ref{lemz32} for any $m \in\mathbb{N}$, for all sufficiently
large $d$,
\[
\ez\bigl[ \bigl( b_d^{k_d^{3/4}} (\mathbf{X}^d_0)
- \ez[ b_d^{k_d^{3/4}} (\mathbf{X}^d_0)]\bigr)^{2m} \bigr] \leq K_m
k_d^{3m/4},
\]
where $K_m = \prod_{j=1}^m (2j-1) +1$. Since $k_d \geq[d^{\beta}]$,
the right-hand side of (\ref{eq32a}) converges to 0 as $\dr$ by
taking $m > 4 \kappa/\beta$, proving (\ref{eq32c}).

Note that
%
%
\begin{equation} \label{eq32b} d^\kappa\pz
(\mathbf{X}_0^d \notin F_d^2) \leq d^\kappa
\sum_{k=[d^\beta]}^{[d^\delta]} \pz\bigl(\mathbf{X}_0^d \notin
\hat{F}_d^2 (k)\bigr).
\end{equation}
The right-hand side of (\ref{eq32b}) converges to 0 as $\dr$ since
(\ref{eq32c}) holds with $\kappa$ replaced by $\kappa+\delta$.\vadjust{\goodbreak}
\end{pf}

Before considering $F_d^3$, the distribution of the position of the
components in the rejection region after $[d^\beta]$ iterations, we
introduce a simple random walk on the hypercube (RWH). The biggest
problem in analyzing the RWM or pseudo-RWM algorithm is the
dependence between the components. However, the dependence is weak
and whether or not a proposed move is accepted is dominated by
whether or not the proposed moves lies inside or outside the
hypercube. Therefore we couple the RWM algorithm to the simpler RWH
algorithm.

For $d \geq1$, define the RWH algorithm as follows. Let
$\mathbf{W}_k^d$ denote the position of the RWH algorithm after $k$
iterations. Then
%
%
\begin{equation} \label{eqrwh1} \mathbf{W}_{k+1}^d =
\cases{
\mathbf{W}_k^d + \sigma_d \mathbf{Z}_{k+1}^d, &\quad
if $\mathbf{W}_k^d + \sigma_d \mathbf{Z}_{k+1}^d \in[0,1]^d$, \cr
\mathbf{W}_k^d, &\quad otherwise.}
\end{equation}
That is, the RWH algorithm simply accepts all
proposed moves which remain inside the hypercube and rejects all
proposed\vspace*{1pt} moves outside the hypercube. Define the pseudo-RWH
algorithm in the obvious fashion with $\hat{\mathbf{W}}_k^d =
(\hat{W}_{k,1}^d, \hat{W}_{k,2}^d, \ldots, \hat{W}_{k,d}^d )$
denoting the position of the pseudo-RWH algorithm at iteration $k$.
Then for $1 \leq i \leq d$, $\hat{W}_{k+1,i}^d = \hat{W}_{k,i}^d +
\sigma_d \hat{Z}_{k+1,i}^d$, where $\hat{Z}_{k+1,i}^d \sim U[ (-
\hat{W}_{k,i}^d /\sigma_d) \vee-1, (\hat{W}_{k,i}^d /\sigma_d)
\wedge1]$.

For our purposes it will suffice to consider the coupling of the
pseudo-RWM and pseudo-RWH algorithms over $[d^\delta]$ iterations
and study how the pseudo-RWH algorithm evolves over $[d^\delta]$
iterations. Note that the RWH algorithm coincides with the RWM
algorithm with a uniform target density over the $d$-dimensional
cube, so in this case the coupling is exact.

The components of the pseudo-RWH algorithm behave independently. For
$x \in(0,1)$, let $\omega_d (x) = \pz(0 < x + \sigma_d Z < 1)$ and
for $\mathbf{x}^d \in(0,1)^d$, let $\Omega_d (\mathbf{x}^d) =
\prod_{j=1}^d \omega_d (x_j)$. Then $\Omega_d (\mathbf{x}^d)$ is the
probability that a proposed move from $\mathbf{x}^d$ is accepted in
the RWH algorithm.
%
%
\begin{lem} \label{lem33}
For any $\alpha< \frac{1}{2}$ and $\mathbf{x}^d \in[0,1]^d$, there
exists a coupling such that
%
%
\begin{equation} \label{eq33a} d^\alpha\pz(\mathbf{X}_1^d \neq
\mathbf{W}_1^d | \mathbf{X}_0^d \equiv\mathbf{W}_0^d =
\mathbf{x}^d) \rightarrow0 \qquad\mbox{as } \dr.
\end{equation}
\end{lem}
\begin{pf}
Let $U \sim U[0,1]$; then we can couple $\mathbf{X}_1^d$ and
$\mathbf{W}_1^d$ using $\mathbf{Z}_1^d$ and~$U$ as follows.
Let
\begin{eqnarray*} 
\mathbf{W}_1^d & = &
\cases{
\mathbf{x}^d
+ \sigma_d \mathbf{Z}_1^d, &\quad if $\mathbf{x}^d + \sigma_d
\mathbf{Z}_1^d \in[0,1]^d$, \cr
\mathbf{x}^d, &\quad otherwise,}
\\
\mathbf{X}_1^d & = & \cases{
\mathbf{x}^d + \sigma_d \mathbf{Z}_1^d, &\quad
if $\mathbf{x}^d + \sigma_d \mathbf{Z}_1^d \in[0,1]^d$\vspace*{2pt}\cr
&\quad and $\displaystyle U \leq1 \wedge
\exp\Biggl(\sum_{j=1}^d \{g
(x_j + \sigma_d Z_{1,j} ) - g(x_j) \} \Biggr)$, \vspace*{2pt}\cr
\mathbf{x}^d, &\quad otherwise.}
\end{eqnarray*}
Therefore, $\mathbf{X}_1^d \neq\mathbf{W}_1^d$ if $\mathbf{x}^d +
\sigma_d \mathbf{Z}_1^d \in[0,1]^d$ and $U > 1 \wedge\exp
(\sum_{j=1}^d \{g
(x_j + \sigma_d Z_{1,j} ) - g(x_j) \} )$.
Thus
%
%
\begin{eqnarray} \label{eq33b}
&& d^\alpha\pz(\mathbf{X}_1^d \neq
\mathbf{W}_1^d | \mathbf{X}_0^d \equiv\mathbf{W}_0^d =
\mathbf{x}^d) \nonumber\\
&&\qquad= d^\alpha\pz\Biggl( \mathbf{x}^d + \sigma_d \mathbf{Z}_1^d
\in[0,1]^d,\nonumber\\
&&\hspace*{59pt} U > 1 \wedge\exp\Biggl(\sum_{j=1}^d \{g (x_j +
\sigma_d Z_{1,j} ) - g(x_j) \} \Biggr) \Biggr)
\nonumber\\[-8pt]\\[-8pt]
&&\qquad= d^\alpha\ez\Biggl[ \prod_{j=1}^d 1_{ \{ 0 < x_j + \sigma_d
Z_{1,j} < 1 \}} \nonumber\\
&&\hspace*{73pt}{}\times\Biggl\{1 - 1 \wedge\exp\Biggl( \sum_{j=1}^d \{g
(x_j + \sigma_d Z_{1,j} ) - g(x_j) \} \Biggr) \Biggr\} \Biggr]
\nonumber\\
&&\qquad\leq d^\alpha\ez\Biggl[ \Biggl| \sum_{j=1}^d \{g (x_j +
\sigma_d Z_{1,j} ) - g(x_j) \} \Biggr| \Biggr],\nonumber
\end{eqnarray}
since for all $y \in\mathbb{R}$, $| 1 - \{1 \wedge\exp(y) \}|
\leq|y|$.

By Taylor's theorem, for $1 \leq j \leq d$, there exists $\xi_j^d$
lying between 0 and $\sigma_d Z_{1,j}$ such that
%
%
\begin{equation} \label{eq33c}
g (x_j + \sigma_d Z_{1,j} ) - g(x_j) = g^\prime(x_j) \sigma_d
Z_{1,j} + \frac{g^{\prime\prime} (x_j + \xi_j^d)}{2} (\sigma_d
Z_{1,j})^2.\hspace*{-18pt}
\end{equation}
Since $g (\cdot)$ is continuously twice
differentiable on $(0,1)$, there exists $K < \infty$ such that
%
%
\begin{equation} \label{eq33d}
\Biggl| \sum_{j=1}^d \{g (x_j + \sigma_d Z_{1,j} ) - g(x_j) \}
\Biggr| \leq\Biggl| \frac{l}{d} \sum_{j=1}^d g^\prime(x_j)
Z_{1,j} \Biggr| + \frac{Kl^2}{2 d}.
\end{equation}

Since the components of $\mathbf{Z}_1^d$ are independent, by
Jensen's inequality, (\ref{eq33d}) and $\ez[(X+c)^2] \leq2 \ez
[X^2] + 2 c^2$, for any random variable $X$ and constant $c$, we
have that
%
%
\begin{eqnarray} \label{eq33e}
&& d^\alpha\pz(\mathbf{X}_1^d \neq
\mathbf{W}_1^d | \mathbf{X}_0^d \equiv\mathbf{W}_0^d =
\mathbf{x}^d) \nonumber\\
&&\qquad\leq d^\alpha\Biggl( 2 \Biggl\{ \frac{l^2}{3 d^2}
\sum_{j=1}^d g^\prime(x_j)^2 + \frac{K^2 l^4}{4 d^2} \Biggr\}
\Biggr)^{{1}/{2}} \\
&&\qquad\rightarrow0 \qquad\mbox{as } \dr,\nonumber
\end{eqnarray}
and
the lemma is proved.
\end{pf}
%
%
\begin{coro}
\label{lem35} Fix $0 < \alpha< \frac{1}{2} - \delta$.

For any $\mathbf{x}^d \in[0,1]^d$, there exists a coupling such
that
%
%
\begin{equation} \label{eq35a} d^\alpha\pz\Biggl(
\bigcup_{j=0}^{[d^\delta]} \{\mathbf{X}_j^d \neq\mathbf{W}_j^d \}|
\mathbf{X}_0^d \equiv\mathbf{W}_0^d = \mathbf{x}^d \Biggr)
\rightarrow0 \qquad\mbox{as } \dr.
\end{equation}

Moreover, if $\mathbf{x}^d \in\tilde{F}_d^1$ and $\alpha+ \delta+
\gamma< \frac{1}{2}$, there exists a coupling such that
%
%
\begin{equation} \label{eq35b} d^\alpha\pz\Biggl(
\bigcup_{j=0}^{[d^\delta]} \{\hat{\mathbf{X}}_j^d \neq
\hat{\mathbf{W}}_j^d \} | \mathbf{X}_0^d \equiv\mathbf{W}_0^d =
\mathbf{x}^d \Biggr) \rightarrow0 \qquad\mbox{as } \dr.
\end{equation}
\end{coro}

For $r \geq0$ and $k=0,1,2,\ldots,$ let
%
%
\begin{equation} \label{eqlam1}
\chi_j^d (x_j;r;k) =
\cases{
1, &\quad if
$\hat{W}_{k,j}^d \in R_d^r$ given that $\hat{W}_{0,j}^d = x_j$, \cr
0, &\quad otherwise.}
\end{equation}
Let $q^d (x;r;k) = \ez[\chi_1^d (x;r;k)]$ and let $\lambda_d
(\mathbf{x}^d; r; k) = \sum_{j=1}^d q^d (x_j;r;k)$. Note that the
movement of the components of the pseudo-RWH algorithm are
independent.
%
%

The next stage in the proof is to show that, if
$\hat{\mathbf{X}}_0^d$ is started in $F_d^3$, then after $k_d$
iterations of the pseudo-RWM algorithm has forgotten its starting
value in terms of the total number and position of the components in
$R_d^l$ (the rejection region). Moreover, the total number and
position of the components in $R_d^l$ after $k_d$ iterations of the
pseudo-RWM algorithm are approximately from the stationary
distribution of $\hat{\mathbf{X}}_\cdot^d$. Before defining and
studying $\{ F_d^3 \}$, we require the following lemma and
associated corollary concerning the distribution of the components
in the rejection region after $k_d$ steps.
%
%
\begin{lem}
\label{lem38a} Let $\{k_d\}$ be any sequence of positive integers
satisfying $[d^\beta] \leq k_d \leq[d^\delta]$.

For any sequence of $\{\mathbf{x}^d \}$ such that $\mathbf{x}^d \in
F_d^1 \cap F_d^2$,
%
%
\begin{equation} \label{eq38aa} d^{2\gamma} \sum_{i=1}^d q^d
(x_i^d; l ; k_d)^2 \rightarrow0 \qquad\mbox{as } \dr.
\end{equation}
Also for all $0 < x < 1$,
%
%
\begin{equation} \label{eq38ab} d^{2 \gamma} q^d
(x; l ; k_d) \rightarrow0 \qquad\mbox{as } \dr.
\end{equation}
\end{lem}
\begin{pf}
Fix $\mathbf{x}^d \in F_d^1 \cap F_d^2$ and set
$\hat{\mathbf{W}}_0^d = \mathbf{x}^d$.

To prove (\ref{eq38aa}) and (\ref{eq38ab}) we couple the components
of $\hat{\mathbf{W}}_t^d$ to a simple reflected random walk process
$\{ S_t^d; t \geq0 \}$. Set $S_0^d = x$ for some $0 < x <1$. Let
$\tilde{Z}_1, \tilde{Z}_2, \ldots$ be i.i.d. according to $U[-1,1]$.
For $t \geq1$, set $S_{t+1}^d = S_t^d + \sigma_d \tilde{Z}_{t+1}$
with reflection at the boundaries 0 and 1 so that $S_t^d \in(0,1)$.
For $x \in(0,1)$, let $p^d (x; l, k_d) = \pz(S_{k_d}^d \in R_d^r |
S_0^d = x)$.\vadjust{\goodbreak}

Consider $\hat{W}_{t,1}^d$ with identical arguments applying for
the other components of~$\hat{\mathbf{W}}_t^d$. Since $k_d
\sigma_d \rightarrow0$ as $\dr$, we assume that $d$ is such that
$(k_d +1) \sigma_d < \frac{1}{2}$. Then
%
%
\begin{equation} \label{eq38ac} \pz\bigl(\hat{W}_{k_d,1}^d \in R_d^r |
\hat{W}_{0,1}^d \in\bigl( ( k_d +1) \sigma_d, 1- (k_d +1) \sigma_d\bigr)\bigr) =0.
\end{equation}
For $x \in(0, (k_d+1) \sigma_d) \cup(1- (k_d+1)
\sigma_d,1)$ we can couple $S_t^d$ and $\hat{W}_{t,1}^d$
such that
%
%
\begin{equation} \label{eq38ad} q^d (x;l;k_d) \leq p^d (x; l ;
k_d).
\end{equation}
For $\sigma_d < y < 1 - \sigma_d$, if $S_t^d = \hat{W}_{t,1}^d = y$,
then set $S_{t+1}^d = \hat{W}_{t+1,1}^d = y + \sigma_d
\tilde{Z}_{t+1}$. Now if $y < \sigma_d$ $(y > 1 - \sigma_d)$,
$\tilde{Z}_{t+1}$ and $\hat{Z}_{t+1,1}^d$ can be coupled such that,
if $S_t^d = \hat{W}_{t,1}^d = y$, then $S_{t+1}^d \leq
\hat{W}_{t+1,1}^d$ $(S_{t+1}^d \geq\hat{W}_{t+1,1}^d)$.
Furthermore, for $y_1 < y_2 < 1/2$ $(y_1 > y_2 > 1/2)$, the above
coupling can be extended to give, if $S_t^d = y_1$ and
$\hat{W}_{t,1}^d = y_2$, then $S_{t+1}^d < W_{t+1}^d$ $(S_{t+1}^d >
W_{t+1}^d)$. Since in $k_d$ iterations either process can move at
most a distance $k_d \sigma_d$, (\ref{eq38ad}) follows from the
above coupling.

Without loss of generality, we assume that $0 < x < (k_d +1)
\sigma_d$ [symmetry arguments apply for $1- (k_d +1) \sigma_d < x <
1$]. By the reflection principle,
%
%
\begin{eqnarray}
\label{eq38ae} p^d (x;l; k_d) &=& \pz\Biggl( - \sigma_d < x +
\sigma_d \sum_{i=1}^{k_d} \tilde{Z}_i < \sigma_d \Biggr)
\nonumber\\[-8pt]\\[-8pt]
&=& \pz\Biggl( - 1 < \frac{x}{\sigma_d} + \sum_{i=1}^{k_d}
\tilde{Z}_i < 1 \Biggr).\nonumber
\end{eqnarray}
By the Berry--Ess\'een theorem, there exists a positive constant,
$K_1 < \infty$ say, such that for all $z \in\mathbb{R}$,
%
%
\begin{equation} \label{eq38af} \Biggl| \pz\Biggl(
\sqrt{\frac{3}{k_d}} \sum_{i=0}^{k_d -1} Z_i \leq z \Biggr) - \Phi
(z) \Biggr| \leq\frac{K_1}{\sqrt{k_d}},
\end{equation}
where
$\Phi(\cdot)$ denotes the c.d.f. of a standard normal. Therefore it
follows from~(\ref{eq38ae}) and (\ref{eq38af}) that there exists a
positive constant, $K_2 < \infty$ say, such that for all $x \in
(0,1)$,
%
%
\begin{equation}
\label{eq38ag} p^d (x;l; k_d) \leq\frac{K_2}{\sqrt{k_d}}.
\end{equation}

By Hoeffding's inequality, for any $\varepsilon> 0$, 
%
%
\begin{eqnarray} \label{eq38ah}
\pz\Biggl( \Biggl| \sum_{i=1}^{k_d} \tilde{Z}_i \Biggr| > \varepsilon
k_d^{3/4} \Biggr) &\leq& 2 \exp\biggl( - \frac{2 ( \varepsilon
k_d^{3/4})^2}{2^2 k_d}
\biggr) \nonumber\\[-8pt]\\[-8pt]
&=& 2 \exp\bigl(- \varepsilon^2 \sqrt{k_d}/2\bigr).\nonumber
\end{eqnarray}
Hence for
$k_d^{3/4}/d< x < (k_d +1)/d$, by taking $\varepsilon= 1/2l$ in
(\ref{eq38ah}), we have that
%
%
\begin{eqnarray} \label{eq38ai} d p^d (x; l ; k_d)
& =& d \pz\Biggl( \Biggl| x + \sigma_d
\sum_{i=1}^{k_d} \tilde{Z}_i \Biggr| < \sigma_d \Biggr) \nonumber
\\
&\leq&
d \pz\Biggl( \Biggl| \sigma_d \sum_{i=1}^{k_d} \tilde{Z}_i
\Biggr|
> \frac{k_d^{3/4}}{2 d} \Biggr) \\
&\leq&
2 d \exp\biggl(- \frac{\sqrt{k_d}}{8 l^2} \biggr) \rightarrow0
\qquad\mbox{as } \dr.\nonumber
\end{eqnarray}
Furthermore, note that for $(k_d +1) \sigma_d < x < 1- (k_d +1)
\sigma_d$, $p(x;l;k_d)=0$.

Then (\ref{eq38ab}) follows immediately from (\ref{eq38ad}) and the
above bounds for $p (x;l;k_d)$ since $d^{2 \gamma} /\sqrt{k_d}
\rightarrow0$ as $\dr$.

Finally, for $\mathbf{x}^d \in F_d^1 \cap F_d^2$, it follows from
(\ref{eq38ad}), (\ref{eq38af}) and (\ref{eq38ag}) that there exists
$K_3 < \infty$ such that
%
%
\begin{eqnarray} \label{eq38j}\quad
d^{2 \gamma} \sum_{i =1}^d q^d
(x_i^d; l
; k_d)^2 & \leq& d^{2 \gamma} \sum_{i =1}^d p^d (x_i^d; l ; k_d)^2
\nonumber\\[-8pt]\\[-8pt]
& \leq& d^{2 \gamma} \biggl\{ K_3 k_d^{3/4} \biggl(
\frac{K_2}{\sqrt{k_d}} \biggr)^2 + 2 d \exp\biggl(-
\frac{\sqrt{k_d}}{8 l^2} \biggr) \biggr\}\nonumber
\end{eqnarray}
with the right-hand side of (\ref{eq38j})
converging to 0 as $\dr$.
\end{pf}
%
%
\begin{coro}
\label{lem37} For any $m \geq2$, any sequence $\{ r_d \}$
satisfying \mbox{$0 \leq r_d \leq l$} and any sequence of positive integers
$\{k_d \}$ satisfying $[d^\beta] \leq k_d \leq[d^\delta]$, there
exists $K < \infty$, such that for all $d \geq1$,
%
%
\begin{equation} \label{eq37a} \ez[ q^d (X_{0,1}^d;r_d; k_d)^m] \leq
K d^{- (1 + \beta m/8)}.
\end{equation}
\end{coro}
\begin{pf}
Fix $m \geq2$. Note that
%
%
\begin{eqnarray} \label{eq37aa} \ez[ q^d (X_{0,1}^d;r_d; k_d)^m] &
\leq&
\ez[ q^d (X_{0,1}^d;l; k_d)^m] \nonumber\\
&=& \int_0^1 q^d (x; l ; k_d)^m f (x) \,dx \nonumber\\[-8pt]\\[-8pt]
&=& \int_{R_d^{k_d^{3/4}}} q^d (x; l ; k_d)^m f (x) \,dx\nonumber\\
&&{} +
\int_{(R_d^{k_d^{3/4}} )^C} q^d (x; l ; k_d)^m f (x)
\,dx.\nonumber
\end{eqnarray}
The two terms on the right-hand side of (\ref{eq37aa})
are bounded using (\ref{eq38ag}) and (\ref{eq38ai}),
respectively.
Thus it follows from the proof of Lemma~\ref{lem38a} that there
exist constants $K_1, K_2 < \infty$ such that, for all $d \geq1$,
%
%
\begin{eqnarray} \label{eq37b}
\ez[ q^d (X_{0,1}^d;r_d; k_d)^m] & \leq& \int_{R_d^{k_d^{3/4}}}
\biggl( \frac{K_1}{\sqrt{k_d}} \biggr)^m f (x) \,dx\nonumber\\
&&{}
+ \int_{(R_d^{k_d^{3/4}} )^C} \biggl\{ 2 \exp\biggl( -
\frac{\sqrt{k_d}}{8l^2} \biggr) \biggr\}^m f (x) \,dx
\nonumber\\
& \leq& \pz(
X_{0,1}^d \in R_d^{k_d^{3/4}} ) \biggl( \frac{K_1}{\sqrt{k_d}}
\biggr)^m\\
&&{} + \pz( X_{0,1}^d \notin R_d^{k_d^{3/4}} )
\times2 \exp\biggl( - \frac{\sqrt{k_d}}{8l^2}
\biggr) \nonumber\\
& \leq& K_2 \frac{k_d^{3/4}}{d} k_d^{-m/2} + 2 \exp\biggl( -
\frac{\sqrt{k_d}}{8l^2} \biggr).\nonumber
\end{eqnarray}
The corollary
follows from (\ref{eq37b}) since $m \geq2$ and $k_d \geq
[d^\beta]$.
\end{pf}

We are now in position to define $\{F_d^3 \}$. For any $0 \leq r
\leq l$ and $k \in\mathbb{Z}^+$, let
%
%
\begin{equation} \label{eqlam2} \hat{F}_d^3 (r;k) = \{ \mathbf
{x}^d\dvtx|
\lambda_d (\mathbf{x}^d; r; k) - \lambda(r) | < d^{-\gamma}/8 \},
\end{equation}
where $\lambda(r) = f^\ast r (1+ r/2l)$. Let
%
%
\begin{equation} \label{eqlam3} F_d^3 = \Bigl\{ \mathbf{x}^d\dvtx\sup
_{[d^\beta] \leq k_d \leq[d^\delta]}
\sup_{0 \leq r \leq l} | \lambda_d (\mathbf{x}^d; r; k_d) - \lambda
(r) | < d^{-\gamma} \Bigr\}.
\end{equation}
We study $\{ \hat{F}_d^3 (r_d, k_d) \}$ as a prelude to analyzing
$\{F_d^3 \}$ where $r_d$ and $k_d$ are defined in Lemma~\ref{lem36}
below.
%
%
\begin{lem}
\label{lem36} For any sequence $\{r_d \}$ satisfying $0 \leq r_d
\leq l$, any sequence of positive integers $\{ k_d \}$ satisfying
$[d^\beta] \leq k_d \leq[d^\delta]$ and $\kappa> 0$,
%
%
\begin{equation} \label{eq36a} d^\kappa\pz\bigl(\mathbf{X}_0^d \notin
\hat{F}_d^3 (r_d,k_d)\bigr) \rightarrow0 \qquad\mbox{as } \dr.
\end{equation}
\end{lem}
\begin{pf}
By the triangle inequality,
%
%
\begin{eqnarray} \label{eq36b}
&&
d^\kappa\pz\bigl(\mathbf{X}_0^d \notin\hat{F}_d^3 (r_d; k_d)\bigr) \nonumber\\
&&\qquad \leq
d^\kappa\pz\bigl( | \lambda_d (\mathbf{X}_0^d; r_d; k_d) - \ez[
\lambda_d (\mathbf{X}_0^d; r_d; k_d)] | > d^{-\gamma}/16\bigr)
\\
&&\qquad\quad{} + d^\kappa\pz\bigl( | \ez[ \lambda_d
(\mathbf{X}_0^d; r_d; k_d)] - \lambda(r_d)| > d^{-\gamma}/16\bigr).\nonumber
\end{eqnarray}
In turn we show that the two terms on the right-hand
side of (\ref{eq36b}) converge to 0 as $\dr$.

By Markov's inequality, we have that for any $m \in\mathbb{N}$,
%
%
\begin{eqnarray}
\label{eq36c} && d^\kappa\pz\bigl( | \lambda_d (\mathbf{X}_0^d;
r_d;k_d) - \ez[ \lambda_d (\mathbf{X}_0^d; r_d; k_d)] | >
d^{-\gamma}/16\bigr)
\nonumber\\
&&\qquad\leq16^{m} d^{\kappa+m \gamma} \ez\Biggl[ \Biggl( \sum_{j =1}^d
\{ q^d (X_{0,j};r_d ;k_d) - \ez[ q^d (X_{0,j};r_d ;k_d)] \}
\Biggr)^m
\Biggr] \nonumber\\[-8pt]\\[-8pt]
&&\qquad= 16^m d^{\kappa+m \gamma} \sum_{i_1 =1}^d \cdots\sum_{i_m =1}^d
\ez\Biggl[ \prod_{j=1}^m \{ q^d (X_{0,i_j};r_d ;k_d)\nonumber\\
&&\hspace*{167pt}{} - \ez[ q^d
(X_{0,i_j};r_d ;k_d)] \} \Biggr].\nonumber
\end{eqnarray}
Since the components of $\mathbf{X}_0^d$ are independent and
identically distributed, we have for any $\{ i_1, i_2, \ldots, i_m
\}$, there exists $1 \leq J \leq m$ and $l_1, l_2, \ldots, l_J \geq
1$ with $l_1 + l_2 + \cdots+ l_J =m$ such that
%
%
\begin{eqnarray}
\label{eq36d}
&&\ez\Biggl[ \prod_{j=1}^m \{ q^d (X_{0,i_j};r_d ;k_d) -
\ez[ q^d (X_{0,i_j};r_d ;k_d)] \} \Biggr]\nonumber\\[-8pt]\\[-8pt]
&&\qquad = \prod_{j=1}^J \ez
\bigl[ \{ q^d (X_{0,1};r_d ;k_d) - \ez[ q^d (X_{0,1};r_d ;k_d)]
\}^{l_j} \bigr].\nonumber
\end{eqnarray}
Note that if any $l_j =1$, then the right-hand side of (\ref{eq36d})
is equal to 0. By Corollary~\ref{lem37}, if $l_1, l_2, \ldots, l_J
\geq2$, there exists $K_1 < \infty$ such that the right-hand side of
(\ref{eq36d}) is less than or equal to $ \prod_{j=1}^J \{ K_1
d^{-(1+l_j\beta/8)} \} = K_1^J d^{-J} d^{-m \beta/8}$. Furthermore,
there exists $K_2 < \infty$ such that for any \mbox{$1 \leq J \leq m$} and
$l_1, l_2, \ldots, l_J \geq2$, there are at most $K_2 d^J$
configurations of $\{i_1, i_2, \ldots, i_m \}$ such that for $j=1,2
,\ldots, J$, $l_j$~of the components are the same. Therefore there
exists $K < \infty$ such that
%
%
\begin{eqnarray} \label{eq36e}
&&\sum_{i_1 =1}^d \cdots\sum_{i_m =1}^d
\ez\Biggl[ \prod_{j=1}^m \{ q^d (X_{0,i_j};r_d ;k_d) - \ez[ q^d
(X_{0,i_j};r_d ;k_d)] \} \Biggr]\nonumber\\[-8pt]\\[-8pt]
&&\qquad \leq K d^{-m \beta/8}.\nonumber
\end{eqnarray}
Taking $m > \kappa/(\beta/8 -\gamma)$, it follows from
(\ref{eq36e}) that the right-hand side of~(\ref{eq36c}) converges to
0 as $\dr$.

The lemma follows by showing that for all sufficiently large $d$,
%
%
\begin{equation} \label{eq36f} |\ez[ \lambda_d (\mathbf{X}_0^d;
r_d; k_d)] - \lambda(r_d)|
\leq d^{-\gamma}/16.
\end{equation}

Note that
%
%
\begin{eqnarray} \label{eq36g}
\ez[ \lambda_d (\mathbf{X}_0^d; r_d; k_d)] &=& d \ez[ q^d
(X_{0,1};r;k_d)] \nonumber\\
&=& d \int_0^{k_d^{3/4}/d} q^d (x;r_d; k_d) f (x)
\,dx\nonumber\\[-8pt]\\[-8pt]
&&{} + d
\int_{k_d^{3/4}/d}^{1-k_d^{3/4}/d} q^d (x;r_d; k_d) f (x) \,dx\nonumber\\
&&{} + d \int_{1-k_d^{3/4}/d}^1 q^d (x;r_d; k_d) f (x)
\,dx.\nonumber
\end{eqnarray}
By (\ref{eq38ai}), the second
integral on the right-hand side of (\ref{eq36g}) is bounded above by
$d \times2 \exp(-\sqrt{k_d}/8l^2) \rightarrow0$ as $\dr$.
Let $f_\star= {\sup_{0 \leq x \leq1}} |f^\prime(x) |$. Then by
Taylor's theorem, for $0 \leq x \leq k_d^{3/4} /d$,
%
%
\begin{equation} \label{eq36ga} | f(x) - f(0)| \leq x \sup_{0 \leq y
\leq x} f^\prime(y) \leq f_\star k_d^{3/4} /d.
\end{equation}
Thus
%
%
\begin{eqnarray} \label{eq36h}\quad
&& \biggl|d \int_0^{k_d^{3/4}/d} q^d (x;r_d; k_d) f (x) \,dx -
f(0) d \int_0^{k_d^{3/4}/d} q^d (x;r_d; k_d) \,dx \biggr|
\nonumber\\[-8pt]\\[-8pt]
&&\qquad\leq d \times f_\star\frac{k_d^{3/4}}{d} \times
\int_0^{k_d^{3/4}/d} q^d (x;r_d; k_d) \,dx.\nonumber
\end{eqnarray}
Similarly, we have that
%
%
\begin{eqnarray} \label{eq36i}\qquad
& &\biggl|d \int_{1-k_d^{3/4}/d}^1 q^d (x;r_d; k_d) f (x) \,dx -
f(1)
d \int_{1-k_d^{3/4}/d}^1 q^d (x;r_d; k_d) \,dx \biggr|
\nonumber\\[-8pt]\\[-8pt]
&&\qquad\leq d \times f_\star\frac{k_d^{3/4}}{d} \times
\int_{1-k_d^{3/4}/d}^1 q^d (x;r_d; k_d) \,dx.\nonumber
\end{eqnarray}
By
symmetry, $q^d (1-x;r_d;k_d) = q^d (x;r_d;k_d)$, so
%
%
\begin{equation} \label{eq36j}
d^\gamma\biggl|\ez[ \lambda_d (\mathbf{X}_0^d; r_d; k_d)] - 2
f^\ast d \int_0^1 q^d (x;r_d;k_d) \,dx \biggr| \rightarrow0
\qquad\mbox{as } \dr.\hspace*{-35pt}
\end{equation}

Since $\int_0^1 \omega_d (y) \,dy \geq1 - 2 \sigma_d$, using Lemma
\ref{lem38a}, (\ref{eq38ab}), we have that, for all sufficiently
large $d$,
%
%
\begin{eqnarray} \label{eq36k}
& & d^\gamma\biggl| d \int_0^1 q^d (x;r_d;k_d) \,dx -
d \int_0^1 q (x;r_d;k_d) \frac{\omega_d (x)}{\int_0^1 \omega_d (y)
\,dy} \,dx \biggr| \nonumber\\
&&\qquad\leq4 d^{1 + \gamma} \int_0^{\sigma_d} q^d (x;r_d;k_d)
\,dx\nonumber\\[-8pt]\\[-8pt]
&&\qquad\quad{} +
d \int_{\sigma_d}^{1-\sigma_d} \biggl\{ \frac{1}{ \int_0^1 \omega_d
(y)
\,dy } - 1 \biggr\} q^d (x;r_d;k_d) \,dx \nonumber\\
&&\qquad\leq4 d^{1 + \gamma} \sigma_d d^{-2 \gamma} + d^{1 + \gamma}
\int_0^1 \frac{2 \sigma_d}{\int_0^1 \omega_d (y) \,dy } q^d
(x;r_d;k_d) \,dx.\nonumber
\end{eqnarray}

Let $p^d (x;r_d;k_d)$ be defined as in Lemma~\ref{lem38a}. Note that
$U[0,1]$ is the stationary distribution of a reflected random walk
on $(0,1)$. Therefore for any $k \geq1$,
%
%
\begin{equation} \label{eq36l} \int_0^1 p^d (x;r_d;k) \,dx = \int
_0^1 p^d (x;r_d;0)
\,dx = \frac{2r_d}{d}.
\end{equation}
Therefore, it follows from
Lemma~\ref{lem38a}, (\ref{eq38ad}), that
%
%
\begin{equation} \label{eq36m} d \int_0^1 q^d
(x;r_d;k_d) \,dx \leq d \int_0^1 p^d (x;r_d;k_d) \,dx = 2 r_d.
\end{equation}
Hence the right-hand side of (\ref{eq36k}) converges
to 0 as $\dr$.

Note that the stationary distribution of a single component of the
pseudo-RWH algorithm has p.d.f. $\omega_d (x) / \int_0^1 \omega_d (y)
\,dy$ $(0 < x <1)$. Therefore
%
%
\begin{eqnarray} \label{eq36n} \qquad
d \int_0^1 q^d (x;r_d;k_d) \frac
{\omega_d (x)}{\int_0^1 \omega_d (y) \,dy} \,dx
&=& d \int_0^1 q^d (x;r_d;0) \frac{\omega_d (x)}{\int_0^1 \omega_d
(y)
\,dy} \,dx \nonumber\\[-8pt]\\[-8pt]
&=& \frac{r_d}{2} \biggl(1 + \frac{r_d}{2l} \biggr) \bigg/ \biggl(1 -
\frac{l}{2d} \biggr).\nonumber
\end{eqnarray}
Finally, combining
(\ref{eq36j}), (\ref{eq36k}) and (\ref{eq36n}), we have that
(\ref{eq36f}) holds and the lemma is proved.
\end{pf}
%
%
\begin{lem}
\label{lema311} For any $\kappa> 0$,
%
%
\begin{equation} \label{eqss2} d^\kappa\pz(\mathbf{X}_0^d
\notin F_d^3) \rightarrow0 \qquad\mbox{as } \dr.
\end{equation}
\end{lem}
\begin{pf}
Fix $\kappa> 0$. Fix a sequence of positive integers
$\{k_d \}$ such that $[d^\beta] \leq k_d \leq[d^\delta]$. Fix $
\theta> \gamma$ and let $\mathcal{S}_d = \{ 0, d^{-\theta}, 2
d^{-\theta}, \ldots, [l d^{\theta}] d^{- \theta},l \}$. Thus the
elements of $\mathcal{S}_d$ are separated by a distance of at most
$d^{-\theta}$.

For any $0 \leq r \leq l$ and $d \geq1$, there exist $\tilde{r}_d,
\hat{r}_d \in\mathcal{S}_d$ such that $\tilde{r}_d \leq r <
\hat{r}_d$ with $\hat{r}_d - \tilde{r}_d \leq d^{-\theta}$. By the
triangle inequality,
%
%
\begin{eqnarray} \label{eqss4}
& & | \lambda_d (\mathbf{X}_0^d; r; k_d) - \lambda(r)| \nonumber\\
&&\qquad\leq\lambda_d (\mathbf{X}_0^d;r;k_d) - \lambda_d
(\mathbf{X}_0^d; \tilde{r}_d;k_d) + | \lambda_d (\mathbf{X}_0^d;
\tilde{r}_d;k_d) - \lambda(\tilde{r}_d)| \nonumber\\
& &\qquad\quad{} + \lambda
(r)- \lambda(\tilde{r}_d) \nonumber\\
&&\qquad\leq\lambda_d (\mathbf{X}_0^d; \hat{r}_d ;k_d) - \lambda_d
(\mathbf{X}_0^d; \tilde{r}_d;k_d) + | \lambda_d (\mathbf{X}_0^d;
\tilde{r}_d;k_d) - \lambda(\tilde{r}_d)| \nonumber\\[-8pt]\\[-8pt]
&&\qquad\quad{} +
\lambda(\hat{r}_d)- \lambda(\tilde{r}_d) \nonumber\\
&&\qquad\leq| \lambda_d (\mathbf{X}_0^d; \hat{r}_d ;k_d) - \ez[
\lambda_d (\mathbf{X}_0^d; \hat{r}_d ;k_d) ]|\nonumber\\
&&\qquad\quad{} + 2 |\lambda_d
(\mathbf{X}_0^d; \tilde{r}_d ;k_d) - \ez[ \lambda_d
(\mathbf{X}_0^d; \tilde{r}_d ;k_d) ] | \nonumber\\
&&\qquad\quad{} + 2 |\lambda(\hat{r}_d)- \lambda(\tilde{r}_d)|.\nonumber
\end{eqnarray}

By Lemma~\ref{lem36}, for any sequence $\{r_d\}$ satisfying $0 \leq
r_d \leq l$,
%
%
\begin{equation} \label{eqss3} \quad
d^{\kappa+ \theta+
\delta} \pz\biggl(| \lambda_d (\mathbf{X}_0^d;r_d ; k_d) - \lambda
(r_d) | > \frac{d^{-\gamma}}{8} \biggr) \rightarrow0 \qquad
\mbox{as } \dr.
\end{equation}
Hence
%
%
\begin{equation} \label{eqss5} \quad
d^{\kappa+
\delta} \pz\biggl( \max_{r \in\mathcal{S}_d} | \lambda_d
(\mathbf{X}_0^d;r ; k_d) - \lambda(r) | > \frac{d^{-\gamma}}{8}
\biggr) \rightarrow0 \qquad\mbox{as } \dr.
\end{equation}
For all sufficiently large $d$,
%
%
\begin{equation} \label{eqss6} \sup_{0 \leq r_d,s_d \leq l,
|r_d-s_d| \leq
d^{-\theta}} | \lambda(r_d) - \lambda(s_d)| \leq
\frac{d^{-\gamma}}{16}.
\end{equation}
Therefore it follows from (\ref{eqss4}), (\ref{eqss5}) and (\ref
{eqss6}) that
%
%
\begin{equation} \label{eqss7} d^{ \kappa+ \delta} \pz\Bigl(
\sup_{0 \leq r \leq l} | \lambda_d
(\mathbf{X}_0^d;r ; k_d) - \lambda(r) | > d^{-\gamma} \Bigr)
\rightarrow0 \qquad\mbox{as } \dr.\hspace*{-20pt}
\end{equation}
Since
(\ref{eqss7}) holds for any sequence $\{ k_d \}$ satisfying
$[d^\beta] \leq k_d \leq[d^\delta]$, the lemma follows since
%
%
\begin{eqnarray} \label{eqss8} && d^{ \kappa} \pz\Bigl( \sup
_{[d^\beta] \leq k \leq[d^\delta]} \sup_{0 \leq r \leq l} | \lambda_d
(\mathbf{X}_0^d;r ; k) - \lambda(r) | > d^{-\gamma} \Bigr)
\nonumber\\[-8pt]\\[-8pt]
&&\qquad\leq d^{ \kappa} \sum_{k=[d^\beta]}^{[d^\delta]}
\pz\Bigl( \sup_{0 \leq r \leq l} | \lambda_d (\mathbf{X}_0^d;r ; k)
- \lambda(r) | > d^{-\gamma} \Bigr).\nonumber
\end{eqnarray}
\upqed\end{pf}

Finally, we consider
%
%
\begin{equation} \label{eqss41} F_d^4 = \Biggl\{ \mathbf{x}^d ;
\Biggl|
\frac{1}{d} \sum_{j=1}^d g^\prime(x_j)^2 - \ez[g^\prime(X_1)^2]
\Biggr| < d^{-{1}/{8}} \Biggr\}.
\end{equation}
The sets
$\{F_d^4 \}$ mirror the sets $\{F_n \}$ in~\cite{RGG} and are used
when considering $G_d^\delta H (\mathbf{x}^d)$ and
$\hat{G}_d^{\delta, \pi} H (\mathbf{x}^d)$ but play no role in
analyzing $P_d$.
%
%
\begin{lem}
\label{lem320} For any $\kappa> 0$,
%
%
\begin{equation} \label{eqss42} d^\kappa\pz(\mathbf{X}^d_0
\notin F_d^4) \rightarrow0 \qquad\mbox{as } \dr.
\end{equation}
\end{lem}
\begin{pf}
Let $g^\ast= \sup_{0 \leq y \leq1} |g^\prime(y)|$
and fix $\kappa> 0$. Then by Hoeffding's inequality,
%
%
\begin{eqnarray} \label{eqss43}
d^\kappa\pz(\mathbf{X}^d_0 \notin F_d^4) &=& d^\kappa\pz\Biggl(
\Biggl| \sum_{j=1}^d g^\prime(X_{0,j})^2 - d \ez[ g^\prime
(X_{0,1})^2] \Biggr| > d^{7/8} \Biggr) \nonumber\\[-8pt]\\[-8pt]
& \leq& d^\kappa\times2 \exp\biggl( - \frac{2 d^{7/4}}{d
(g^\ast)^4} \biggr) \rightarrow0 \qquad\mbox{as } \dr.\nonumber
\end{eqnarray}
\upqed\end{pf}

Finally we are in position to consider $\{F_d \}$ and $\{\tilde{F}_d
\}$. Recall that, for $d \geq1$, $F_d = F_d^1 \cap F_d^2 \cap
F_d^3 \cap F_d^4$ and
\[
\tilde{F}_d = \Biggl\{ \mathbf{x}^d; \pz
\biggl( \bigcup_{j=0}^{[d^\delta]} \hat{\mathbf{X}}_j^d \notin
F_d |
\hat{\mathbf{X}}_0^d = \mathbf{x}^d \biggr) \leq d^{-3} \Biggr\}.
\]

Combining Lemmas~\ref{lem31},~\ref{lem32},~\ref{lema311}
and~\ref{lem320}, we have the following theorem.
%
%
\begin{theorem}
\label{lem321} For any $\kappa> 0$,
%
%
\begin{equation} \label{eqss45} d^\kappa\pz(\mathbf{X}_0^d \notin F_d)
\rightarrow
0 \qquad\mbox{as } \dr.
\end{equation}
Hence, by Lemma
\ref{lem34}, for any $\kappa> 0$,
%
%
\begin{equation} \label{eqss46} d^\kappa\pz(\mathbf{X}_0^d \notin\tilde
{F}_d) \rightarrow
0 \qquad\mbox{as } \dr.
\end{equation}

Also using the couplings outlined above, we have that
%
%
\begin{equation} \label{eqss47} \pz\Biggl( \bigcup
_{j=0}^{[d^\delta]}
\{ \hat{\mathbf{W}}_j^d \notin F_d \} |\hat{\mathbf{W}}_0^d \in
\tilde{F}_d \Biggr) \rightarrow0 \qquad\mbox{as } \dr.
\end{equation}
\end{theorem}
%
%

\section{\texorpdfstring{Proof of $P_d|\mathbf{X}_0^d =\mathbf{x}^d \convp\exp(-l f^\ast/2)$}
{Proof of Pd|X0d =xd ->p exp(-l f*/2)}} \label{SecPd}

We show that for any sequence $\{\mathbf{x}^d \}$ such that
$\mathbf{x}^d \in\tilde{F}_d$,
%
%
\begin{equation} \label{eqspd1}
P_d | \mathbf{X}_0^d =\mathbf{x}^d \convp\exp(-l f^\ast/2)
\qquad\mbox{as } \dr.
\end{equation}
The key result is
Lemma~\ref{lem311} which states that after $k_d$ iterations, the
configuration of the components in the rejection region $R_d^l$
resemble the configuration of the points of a Poisson point process
with rate $\lambda(r) = f^\ast r (1+ r/2l)$ on the interval $[0,l]$.

For any $n \in\mathbb{N}$ and $1 \leq i \leq n$, let
\[
S_n^d (\mathbf{x}^d; i ; k) = \sum_{j=1}^d \bigl\{ \chi_i^d
(x_j; il/n ;k) - \chi_j^d \bigl(x_j; (i-1)l/n ;k\bigr) \bigr\}
\]
with
\[
\mathbf{S}_n^d (\mathbf{x}^d; k) = (S_n^d (\mathbf{x}^d; 1 ;
k), S_n^d (\mathbf{x}^d; 2 ; k), \ldots, S_n^d (\mathbf{x}^d; n ;
k)).
\]
Let $\mathbf{S}_n = (S_n^1, S_n^2, \ldots, S_n^n)$ where the
components of $\mathbf{S}_n$ are independent Poisson random
variables with $S_n^i \sim\operatorname{Po} (\lambda_{n,i})$ and
\[
\lambda_{n,i} = \lambda(il/n) - \lambda\bigl((i-1)l/n\bigr)
\qquad(1 \leq i \leq n).
\]
%

\begin{lem}
\label{lem311} For any $n \in\mathbb{N}$, any sequence of positive
integers $\{ k_d \}$ satisfying $[d^\beta] \leq k_d \leq[d^\delta]$
and $\mathbf{x}^d \in F_d$,
\[
\mathbf{S}_n^d
(\mathbf{x}^d; k_d) \convd\mathbf{S}_n \qquad\mbox{as }
\dr.
\]
\end{lem}
\begin{pf}
Fix $n \in\mathbb{N}$ and $\mathbf{x}^d \in F_d$. Let
\[
\check{\mathbf{S}}_n^d
(\mathbf{x}^d; k_d) = ( \check{S}_n^d (\mathbf{x}^d; 1 ; k_d),
\check{S}_n^d (\mathbf{x}^d; 2 ; k_d), \ldots, \check{S}_n^d
(\mathbf{x}^d; n ; k_d)),
\]
where for $1 \leq i \leq n$,
$\check{S}_n^d (\mathbf{x}^d; i ; k_d)$ are independent Poisson
random variables with means
\[
\lambda_{n,i}^d (\mathbf{x}^d; k_d) =\lambda_d (\mathbf{x}^d; il/n;
k_d) - \lambda_d \bigl(\mathbf{x}^d; (i-1)l/n; k_d\bigr).
\]

The lemma is proved by showing that
%
%
\begin{eqnarray} \label{eq311b} d_{\mathit{TV}} ( \tilde{\mathbf{S}}_n^d
(\mathbf{x}^d; k_d) , \mathbf{S}_n) &\leq& d_{\mathit{TV}} ( \mathbf{S}_n^d
(\mathbf{x}^d; k_d) , \check{\mathbf{S}}_n^d (\mathbf{x}^d; k_d))\nonumber\\
&&{} +
d_{\mathit{TV}} ( \check{\mathbf{S}}_n^d (\mathbf{x}^d; k_d) , \mathbf{S}_n)
\\
&\rightarrow& 0
\qquad\mbox{as } \dr.\nonumber
\end{eqnarray}

By~\cite{barb}, Theorem 1,
%
%
\begin{equation} \label{eq311c} d_{\mathit{TV}} (
\tilde{\mathbf{S}}_n^d (\mathbf{x}^d; k_d) , \check{\mathbf{S}}_n^d
(\mathbf{x}^d; k_d) ) \leq\sum_{i=1}^d q^d (x_i; l ; k_d)^2.
\end{equation}
By Lemma~\ref{lem38a}, (\ref{eq38aa}) the right-hand side of
(\ref{eq311c}) converges to 0 as $\dr$.

For the second term on the right-hand side of (\ref{eq311b}), it
suffices to show that
\[
\check{\mathbf{S}}_n^d (\mathbf{x}^d;
k_d) \convd\mathbf{S}_n \qquad\mbox{as } \dr.
\]
(For discrete random variables convergence in
distribution and convergence in total variation distance are
equivalent; see~\cite{BHJ}, page 254.)

The components of $ \check{\mathbf{S}}_n^d (\mathbf{x}^d; k_d)$ and
$\mathbf{S}_n$ are independent, and therefore it is sufficient to
show that, for all $1 \leq i \leq n$,
%
%
\begin{equation} \label{eq311d} \check{S}_n^d (\mathbf{x}^d;
i; k_d) \convd S_{n,i} \qquad\mbox{as } \dr.
\end{equation}
For all $1 \leq i \leq n$, (\ref{eq311d}) holds, if
%
%
\begin{equation} \label{eq311e} \lambda_{n,i}^d (\mathbf{x}^d; k_d)
\rightarrow\lambda_{n,i} \qquad\mbox{as } \dr.
\end{equation}
Therefore the lemma follows from (\ref{eq311e}) since $[d^\beta]
\leq k_d \leq[d^\delta]$ and $\mathbf{x}^d \in F_d^3$. [See
(\ref{eqn11f3}) for the construction of $\{F_d^3 \}$.]\vadjust{\goodbreak}
\end{pf}

Lemma~\ref{lem311} is the key result stating that if the pseudo-RWH
process is started from the set $F_d$, then after $[d^\beta]$
iterations the distribution of the components in the rejection
region are approximately given by $\mathbf{S}_n$. We show that
studying the pseudo-RWH algorithm over $[d^\delta]$ iterations
suffices in analyzing $T_d (\pi) = \frac{1}{[d^\delta]}
\sum_{j=0}^{[\pi d^\delta-1]} M_j (J_d (\hat{\mathbf{X}}_j^d))$.
Note that $P_d$ satisfies
%
%
\begin{equation} \label{eqtt1} T_d (P_d) \leq1 < T_d (P_d +
1/[d^\delta]).
\end{equation}

Let $\hat{T}_d (\pi) =\frac{1}{[d^\delta]} \sum_{j=0}^{[\pi
d^\delta
-1]} M_j (\Omega_d (\hat{\mathbf{W}}_j^d))$. Before establishing a
coupling between $T_d (\pi)$ and $\hat{T}_d (\pi)$, we give a simple
coupling for geometric random variables.
%
%
\begin{lem}
\label{lem310} Suppose that $0 \leq q < p \leq1$ and that $X$ and
$Y$ are independent geometric random variables with success
probabilities $p$ and $q$, respectively, that is, $X \sim M (p)$
and $Y \sim M (q)$. Let $A$ be a Bernoulli random variable with $\pz
(A =1) =q/p$ and $Z \sim M (q)$. Then if $A$, $X$, $Y$ and $Z$ are
mutually independent,
%
%
\begin{equation} \label{eq310a}
Y \eqd X + (1-A) Z.
\end{equation}
Therefore there exists a coupling of $X$ and $Y$ such that
%
%
\begin{equation} \label{eq310b} \pz(X \neq Y) = P (A=0) =
\frac{p-q}{p}.
\end{equation}
\end{lem}
%
%
\begin{lem}
\label{lem312} For any $0 < \pi\leq1$ and $\mathbf{x}^d \in
\tilde{F}_d$, there exists a coupling of $T_d (\pi)$ and $\hat{T}_d
(\pi)$ such that
%
%
\begin{equation} \label{eq312a}
\pz\bigl( T_d (\pi) \neq\hat{T}_d (\pi) | \hat{\mathbf{X}}_0^d \equiv
\hat{\mathbf{W}}_0^d = \mathbf{x}^d\bigr) \rightarrow0 \qquad
\mbox{as } \dr.
\end{equation}
\end{lem}
\begin{pf}
For $\mathbf{x}^d \in\tilde{F}_d$, by Corollary
\ref{lem35}, we have that
%
%
\begin{equation} \label{eq312b}
\pz\Biggl( \bigcup_{j=0}^{[d^\delta]} \{ \hat{\mathbf{X}}_j^d \neq
\hat{\mathbf{W}}_j^d \} | \hat{\mathbf{X}}_0^d \equiv
\hat{\mathbf{W}}_0^d = \mathbf{x}^d \Biggr) \rightarrow0
\qquad\mbox{as } \dr.
\end{equation}

Suppose that for $j=0,1,\ldots, [d^\delta]$, $\hat{\mathbf{W}}_j^d =
\hat{\mathbf{X}}_j^d \in F_d^1$. Then using Lemma~\ref{lem310},
(\ref{eq310b}), $M_j (J_d (\hat{\mathbf{X}}_j^d))$ and $M_j
(\Omega_d (\hat{\mathbf{W}}_j^d))$ can be coupled such that
%
%
\begin{equation} \label{eq312d}\hspace*{32pt}
\pz\bigl(M_j (J_d (\hat{\mathbf{X}}_j^d)) \neq M_j (\Omega_d
(\hat{\mathbf{W}}_j^d))| \hat{\mathbf{W}}_j^d = \hat{\mathbf{X}}_j^d
\in F_d^1\bigr) \leq\frac{\Omega_d (\hat{\mathbf{X}}_j^d) - J_d
(\hat{\mathbf{X}}_j^d)}{\Omega_d (\hat{\mathbf{X}}_j^d)}.
\end{equation}
Since $\hat{\mathbf{X}}_j^d \in F_d^1$, $\Omega_d
(\hat{\mathbf{X}}_j^d) \geq2^{-\gamma\log d} \geq d^{-\gamma}$,
the right-hand side of (\ref{eq312d}) is less than $d^\gamma\{\Omega_d
(\hat{\mathbf{X}}_j^d) - J_d (\hat{\mathbf{X}}_j^d) \}$. Note that
\[
\pz(\hat{\mathbf{W}}_{j+1}^d \neq\hat{\mathbf{X}}_{j+1}^d|
\hat{\mathbf{W}}_j^d = \hat{\mathbf{X}}_j^d
\in F_d^1) = \Omega_d (\hat{\mathbf{X}}_j^d) - J_d
(\hat{\mathbf{X}}_j^d),
\]
so by Lemma~\ref{lem33} for any $\alpha<
\frac{1}{2}$, $d^{\alpha- \gamma}$ times the right-hand side of
(\ref{eq312d}) converges to 0 as $\dr$. Taking $\alpha$ such that
$\delta+ \gamma< \alpha< \frac{1}{2}$,
%
%
\begin{eqnarray} \label{eq312e}
&&
\sum_{j=0}^{[d^\delta]} \pz\bigl(M_j (J_d (\hat{\mathbf{X}}_j^d)) \neq
M_j (\Omega_d (\hat{\mathbf{W}}_j^d))| \hat{\mathbf{W}}_j^d =
\hat{\mathbf{X}}_j^d \in F_d^1\bigr)\nonumber\\[-8pt]\\[-8pt]
&&\qquad\rightarrow0 \qquad
\mbox{as } \dr.\nonumber
\end{eqnarray}
The lemma then follows from (\ref{eq312b}) and (\ref{eq312e}).
\end{pf}

We show that it suffices to study $\tilde{T}_d (\pi) =
\frac{1}{[d^\delta]} \sum_{j=0}^{[\pi d^\delta-1]} \Omega_d
(\hat{\mathbf{W}}_j^d)^{-1}$. In other words, replace the mean of
the geometric random variables $\{M (\Omega_d
(\hat{\mathbf{W}}_0^d))$, $M (\Omega_d (\hat{\mathbf{W}}_1^d)),
\ldots,M (\Omega_d (\hat{\mathbf{W}}_{[\pi d^\delta-1]}^d)) \}$ by
the mean of the means of the geometric random variables.
%
%
\begin{lem}
\label{lem313} For any $0 < \pi\leq1$ and for any sequence of $\{
\mathbf{x}^d \}$ such that \mbox{$\mathbf{x}^d \in\tilde{F}_d$}, $
\hat{T}_d (\pi) | \hat{\mathbf{W}}_0^d = \mathbf{x}^d \convd\pi
\exp(f^\ast l/2)$ if $\tilde{T}_d (\pi) | \hat{\mathbf{W}}_0^d =
\mathbf{x}^d \convd\pi\exp(f^\ast l/2)$ as $\dr$.
\end{lem}
\begin{pf}
Let $A_d = \bigcup_{j=0}^{[d^\delta]}
\{\hat{\mathbf{W}}_j^d \notin F_d \}$. Then for any $\mathbf{x}^d
\in\tilde{F}_d$, $\pz(A_d | \hat{\mathbf{W}}_0^d = \mathbf{x}^d)
\rightarrow0$ as $\dr$.

For any $\tau\in\mathbb{R}$ with $i = \sqrt{-1}$, the
characteristic function of $\hat{T}_d (\pi)$ conditional upon
$A_d^C$ and $\hat{\mathbf{W}}_0^d = \mathbf{x}^d$ is given by
%
%
\begin{eqnarray} \label{eq313a}\qquad
&& \ez[ \exp(i \tau\hat{T}_d (\pi))
| A_d^C , \hat{\mathbf{W}}_0^d = \mathbf{x}^d] \nonumber\\
&&\qquad= \ez\Biggl[ \prod_{j=0}^{[\pi d^\delta-1]} \ez\biggl[
\exp\biggl( \frac{i \tau}{[d^\delta]} M_j (\Omega_d
(\hat{\mathbf{W}}_j^d)) \biggr) \Big| A_d^C, \{ \hat{\mathbf{W}}^d
\} \biggr] \Big| A_d^C, \hat{\mathbf{W}}_0^d = \mathbf{x}^d
\Biggr] \\
&&\qquad= \ez\Biggl[ \prod_{j=0}^{[\pi d^\delta-1]} \frac{\exp(i
\tau/[d^\delta]) \Omega_d (\hat{\mathbf{W}}_j^d)}{1 - (1 - \Omega_d
(\hat{\mathbf{W}}_j^d)) \exp(i \tau/[d^\delta])} \Big| A_d^C,
\hat{\mathbf{W}}_0^d = \mathbf{x}^d \Biggr].\nonumber
\end{eqnarray}
Conditional upon $A_d^C$, $\Omega_d (\hat{\mathbf{W}}_j^d)^{-1} \leq
2^{\gamma\log d} \leq d^\gamma$. Hence, for all \mbox{$0 \leq j \leq
[\pi d^\delta-1 ]$},
%
%
\begin{equation} \label{eq313b}\qquad
\frac{\exp(i \tau/[d^\delta]) \Omega_d (\hat{\mathbf{W}}_j^d)}{1 -
(1 - \Omega_d (\hat{\mathbf{W}}_j^d)) \exp(i \tau/[d^\delta])} = 1
+ \frac{i \tau}{[d^\delta]} \Omega_d (\hat{\mathbf{W}}_j^d)^{-1}
+ o
\biggl( \frac{1}{[d^\delta]} \biggr).
\end{equation}
Thus $\ez[ \exp
(i \tau\hat{T}_d (\pi)) | A_d^C , \hat{\mathbf{W}}_0^d =
\mathbf{x}^d] $ has the same limit as $\dr$ (should one exist) as
%
%
\begin{equation} \label{eq313c}
\ez\Biggl[ \prod_{j=0}^{[\pi d^\delta-1]} \biggl( 1 + \frac{i
\tau}{[d^\delta]} \Omega_d (\hat{\mathbf{W}}_j^d)^{-1} \biggr)
\Big| A_d^C , \hat{\mathbf{W}}_0^d = \mathbf{x}^d \Biggr],
\end{equation}
which in turn has the same limit as $\dr$ as
%
%
\begin{eqnarray} \label{eq313d}
&&\ez\Biggl[\prod_{j=0}^{[\pi d^\delta-1]} \exp
\biggl( \frac{i \tau}{[d^\delta]} \Omega_d
(\hat{\mathbf{W}}_j^d)^{-1} \biggr) \Big| A_d^C ,
\hat{\mathbf{W}}_0^d = \mathbf{x}^d \Biggr]\nonumber\\[-8pt]\\[-8pt]
&&\qquad = \ez[\exp(i \tau
\tilde{T}_d (\pi)) | A_d^C , \hat{\mathbf{W}}_0^d =
\mathbf{x}^d].\nonumber
\end{eqnarray}
The lemma follows since $\pz( A_d^C |
\hat{\mathbf{W}}_0^d = \mathbf{x}^d) \rightarrow1$ as $\dr$.
\end{pf}

We shall show that $\tilde{T}_d (\pi) \convp\exp(l f^\ast/2)$ as
$\dr$ using Chebyshev's inequality in Lemma~\ref{lem318}. We
require preliminary results concerning\break $\operatorname{cov} (\Omega_d
(\hat{\mathbf{W}}_j^d)^{-1}$, $\Omega_d
(\hat{\mathbf{W}}_{j+k}^d)^{-1} | \hat{\mathbf{W}}_0^d =
\mathbf{x}^d)$ with the key results given in Lem\-ma~\ref{lem317}.
First, however, we introduce useful upper and lower bounds for
$\Omega_d (\mathbf{x}^d)^{-1}$ which allow us to exploit Lemma
\ref{lem311} and prove uniform integrability $\{ \tilde{T}_d (\pi)
\}$.\vspace*{1pt}

For $n \in\mathbb{N}$, $1 \leq i \leq n$ and $\mathbf{x}^d \in
(0,1)^d$, let $\tilde{b}_d^{n,i} (\mathbf{x}^d) = b_d^{il/n}
(\mathbf{x}^d) - b_d^{(i-1)l/n} (\mathbf{x}^d)$ with
$\tilde{\mathbf{b}}_d^n (\mathbf{x}^d) = (\tilde{b}_d^{n,1}
(\mathbf{x}^d), \tilde{b}_d^{n,2} (\mathbf{x}^d), \ldots,
\tilde{b}_d^{n,n} (\mathbf{x}^d))$. For $n \in\mathbb{N}$ and
$\mathbf{s} = (s_1, s_2, \ldots,\break s_n) \in\mathbb{R}^n$, let
%
%
\begin{eqnarray} \label{eqnu1} \check{\nu}_n (\mathbf{s}) &=&
\prod_{j=1}^n \biggl( \frac{1}{2} + \frac{j-1}{2n} \biggr)^{-s_j},
\\
\label{eqnu2} \hat{\nu}_n (\mathbf{s}) &=& \prod_{j=1}^n \biggl(
\frac{1}{2} + \frac{j}{2n} \biggr)^{-s_j}.
\end{eqnarray}
Then for
all $\mathbf{x}^d \in(0,1)^d$,
%
%
\begin{equation} \label{eqnu3} \hat{\nu}_n ( \tilde{\mathbf{b}}_d^n
(\mathbf{x}^d)) \leq\Omega_d (\mathbf{x}^d)^{-1} \leq\check{\nu}_n
( \tilde{\mathbf{b}}_d^n (\mathbf{x}^d)) \leq2^{b_d^l
(\mathbf{x}^d)}.
\end{equation}

%
\begin{lem}
\label{lem314} For any $m \in\mathbb{N}$, any sequence of
$\{\mathbf{x}^d \}$ such that $\mathbf{x}^d \in\tilde{F}_d$ and any
sequence of positive integers $\{ k_d \}$ satisfying $[d^\beta] \leq
k_d \leq[d^\delta]$,
%
%
\begin{equation} \label{eq314a} \ez\bigl[ \bigl( 2^{b_d^l
(\hat{\mathbf{W}}^d_{k_d})} \bigr)^m | \hat{\mathbf
{W}}_0^d =
\mathbf{x}^d \bigr] \rightarrow\exp\bigl( (2^m -1) \lambda(l)
\bigr) \qquad\mbox{as } \dr.\hspace*{-20pt}
\end{equation}
\end{lem}
\begin{pf}
Note that $\{b_d^l (\hat{\mathbf{W}}^d_{k_d}) |
\hat{\mathbf{W}}_0^d = \mathbf{x}^d \} = \sum_{j=1}^d \chi_j
(x_j;l;k_d)$. Then since the $\{ \chi_j (x_j;l;k_d) \}$ are
independent Bernoulli random variables,
%
%
\begin{eqnarray} \label{eq314b}\qquad
\ez\bigl[ \bigl( 2^{b_d^l (\hat{\mathbf{W}}^d_{k_d})}
\bigr)^m | \hat{\mathbf{W}}_0^d = \mathbf{x}^d \bigr] &=&
\prod_{j=1}^d \ez\bigl[ (2^m)^{\chi_j (x_j;l;k_d)} |
\hat{\mathbf{W}}_0^d = \mathbf{x}^d \bigr] \nonumber\\[-8pt]\\[-8pt]
&= & \prod_{j=1}^d \bigl\{ \bigl(1 - q^d (x_j;l;k_d)\bigr) + 2^m q^d
(x_j;l;k_d) \bigr\}.\nonumber
\end{eqnarray}
By Lemma~\ref{lem38a},
(\ref{eq38aa}), for $\mathbf{x}^d \in\tilde{F}_d$, $\sum_{j=1}^d
q^d (x_j;l;k_d)^2 \rightarrow0$ as $\dr$, so the right-hand side of
(\ref{eq314b}) has the same limit as $\dr$ as
%
%
\begin{equation} \label{eq314c} \prod_{j=1}^d \exp\bigl( (2^m-1) q^d
(x_j;l;k_d) \bigr) = \exp\bigl( (2^m -1) \lambda_d
(\mathbf{x}^d;l;k_d) \bigr).
\end{equation}
The lemma follows since
for any $\mathbf{x}^d \in\tilde{F}_d$, $\lambda_d
(\mathbf{x}^d;l;k_d) \rightarrow\lambda(l)$ as $\dr$.
\end{pf}
%
%
\begin{lem}
\label{lem315}Fix $m, n \in\mathbb{N}$. For any sequence $\{
\mathbf{x}^d \}$ such that $\mathbf{x}^d \in F_d$, and any sequence
of positive integers $\{ k_d \}$ satisfying $[d^\beta] \leq k_d \leq
[d^\delta]$, we have that
\begin{eqnarray*}
\ez[ \check{\nu}_n (\tilde{\mathbf{S}}_n^d
(\mathbf{x}^d; k_d))^m ] &\rightarrow&\ez[ \check{\nu}_n
(\mathbf{S}_n)^m] \qquad\mbox{as } \dr, \\
\ez[ \hat{\nu}_n (\tilde{\mathbf{S}}_n^d
(\mathbf{x}^d; k_d))^m] &\rightarrow&\ez[ \hat{\nu}_n^m
(\mathbf{S}_n)^m] \qquad\mbox{as } \dr.
\end{eqnarray*}
\end{lem}
\begin{pf}
By~\cite{Bill79}, Theorem 29.2, and Lemma~\ref{lem311}
%
%
\begin{eqnarray} \label{eq315a} \check{\nu}_n(\tilde{\mathbf
{S}}_n^d (\mathbf{x}^d;
k_d))^m
&\convd&\check{\nu}_n(\mathbf{S}_n)^m \qquad\mbox{as } \dr
, \\
\label{eq315b} \hat{\nu}_n (\tilde{\mathbf{S}}_n^d (\mathbf{x}^d;
k_d))^m &\convd&\hat{\nu}_n (\mathbf{S}_n)^m \qquad\mbox{as
} \dr.
\end{eqnarray}
The lemma follows since (\ref{eqnu3}) and Lemma
\ref{lem314} ensure the uniform integrability of the left-hand sides
of (\ref{eq315a}) and (\ref{eq315b}).
\end{pf}
%
%
\begin{lem}
\label{lem316} For any sequence $\{ \mathbf{x}^d \}$ such that
$\mathbf{x}^d \in F_d$ and sequence of positive integers $\{k_d\}$
satisfying $[d^\beta] \leq k_d \leq[d^\delta]$,
%
%
\begin{equation} \label{eq316a} \ez[ \Omega_d
(\hat{\mathbf{W}}_{k_d}^d)^{-1} | \hat{\mathbf{W}}_0^d =
\mathbf{x}^d] \rightarrow\exp(f^\ast l/2)
\qquad\mbox{as
} \dr.
\end{equation}

For any $\mathbf{x}^d \in\tilde{F}_d$ and sequences of positive
integers $\{i_d \}$ and $\{k_d\}$ satisfying $[d^\beta] \leq k_d
\leq[d^\delta]$ and $i_d + k_d \leq[d^\delta]$,
%
%
\begin{equation} \label{eq316b} \qquad
\ez[ \Omega_d (\hat{\mathbf{W}}_{i_d + k_d}^d)^{-1} |
\hat{\mathbf{W}}_{i_d}^d, \hat{\mathbf{W}}_0^d = \mathbf{x}^d]
\convp\exp(f^\ast l/2) \qquad\mbox{as } \dr.
\end{equation}
\end{lem}
\begin{pf}
An immediate consequence of Lemma~\ref{lem315} is that
\[
\lim_{\dr} \ez[ \check{\nu}_n(\tilde{\mathbf{S}}_n^d (\mathbf
{x}^d; k_d))],
\lim_{\dr} \ez[ \hat{\nu}_n(\tilde{\mathbf{S}}_n^d (\mathbf{x}^d;
k_d))] \rightarrow\exp(f^\ast l/2) \qquad\mbox{as } \nr,
\]
from which (\ref{eq316a}) follows by (\ref{eqnu3}).

By Theorem~\ref{lem321}, (\ref{eqss47}), $\pz(
\hat{\mathbf{W}}_{i_d}^d \in F_d | \hat{\mathbf{W}}_0^d \in
\tilde{F}_d) \rightarrow1$ as $\dr$, so (\ref{eq316b}) follows from
(\ref{eq316a}).
\end{pf}
%
%
\begin{lem}
\label{lem317} For any sequence $\{ \mathbf{x}^d \}$ such that
$\mathbf{x}^d \in\tilde{F}_d$ and any sequences of positive
integers $\{ i_d \}$ and $\{k_d\}$ satisfying $[d^\beta] \leq i_d,
k_d \leq[d^\delta]$,
%
%
\begin{equation} \label{eq317a}\quad
\operatorname{cov} \bigl(\Omega_d (\hat{\mathbf
{W}}_{i_d}^d)^{-1}
, \Omega_d (\hat{\mathbf{W}}_{i_d +k_d}^d)^{-1} |
\hat{\mathbf{W}}_0^d = \mathbf{x}^d \bigr) \rightarrow0
\qquad\mbox{as } \dr
\end{equation}
and
%
%
\begin{eqnarray} \label{eq317b}
&&\operatorname{var} \bigl(\Omega_d
(\hat{\mathbf{W}}_{k_d}^d)^{-1}
| \hat{\mathbf{W}}_0^d = \mathbf{x}^d \bigr) \nonumber\\[-9pt]\\[-9pt]
&&\qquad\rightarrow\exp
\bigl(f^\ast l \bigl\{ 4 \log2 - \tfrac{3}{2} \bigr\} \bigr) -
\exp(f^\ast l ) \qquad\mbox{as } \dr.\nonumber
\end{eqnarray}
\end{lem}
\begin{pf}
Using (\ref{eqnu3}), Lemma~\ref{lem314} and Markov's
inequality, it is straightforward to show that for any $\delta> 0$,
there exists $K < \infty$ such that
%
%
\begin{equation} \label{eq317c}\qquad
\pz\bigl(\Omega_d (\hat{\mathbf{W}}_{j_d}^d)^{-1} > K|
\hat{\mathbf{W}}_0^d = \mathbf{x}^d \bigr) \leq\pz\bigl(2^{b_d^l
(\hat{\mathbf{W}}_{j_d}^d)} > K| \hat{\mathbf{W}}_0^d = \mathbf{x}^d
\bigr) \leq\delta.
\end{equation}
Therefore it follows from Lemma
\ref{lem316} that, for any sequence $\{\mathbf{x}^d \}$ such that
$\mathbf{x}^d \in\tilde{F}_d$,
%
%
\begin{eqnarray} \label{eq317d}
&& \Omega_d
(\hat{\mathbf{W}}_{j_d}^d)^{-1} \{ \ez[ \Omega_d
(\hat{\mathbf{W}}_{j_d +k_d}^d)^{-1} | \hat{\mathbf{W}}_{j_d}^d,
\hat{\mathbf{W}}_0^d = \mathbf{x}^d ] \nonumber\\[-2pt]
&&\hspace*{47pt}\qquad{}- \ez[ \Omega_d
(\hat{\mathbf{W}}_{j_d +k_d}^d)^{-1} | \hat{\mathbf{W}}_0^d =
\mathbf{x}^d ] \} | \hat{\mathbf{W}}_0^d = \mathbf{x}^d
\\[-2pt]
&&\qquad \convp0\qquad \mbox{as } \dr.\nonumber
\end{eqnarray}
The uniform integrability of the left-hand side of
(\ref{eq317d}) follows from (\ref{eqnu3}) and Lemma~\ref{lem314}.
Hence (\ref{eq317a}) follows.

It is straightforward to show that $\ez[\check{\nu}_n
(\mathbf{S}_n)^2], \ez[ \hat{\nu}_n (\mathbf{S}_n)^2] \rightarrow
\exp(f^\ast l (4 \log2 - 3/2))$ as $\nr$. Therefore from
(\ref{eqnu3}) and Lemma~\ref{lem314}, we have that
%
%
\begin{equation} \label{eq317e} \ez[ \Omega_d
(\hat{\mathbf{W}}_{k_d}^d)^{-2} | \hat{\mathbf{W}}_{0}^d =
\mathbf{x}^d] \rightarrow\exp\bigl(f^\ast l (4 \log2 - 3/2)\bigr)
\qquad\mbox{as } \dr.\hspace*{-35pt}
\end{equation}
Then (\ref{eq317b}) follows immediately.\vspace*{-3pt}
\end{pf}

We are now in position to prove Lemma~\ref{lem318}, which is the
final step in proving that for any sequence $\{ \mathbf{x}^d \}$
such that $\mathbf{x}^d \in\tilde{F}_d$, $ P_d | \mathbf{X}_0^d =
\mathbf{x}^d \convp\exp(-f^\ast l/2)$ as $\dr$.\vspace*{-3pt}
%
%
\begin{lem}
\label{lem318} For any $0 < \pi\leq1$ and any sequence $\{
\mathbf{x}^d \}$ such that \mbox{$\mathbf{x}^d \in\tilde{F}_d$},
%
%
\begin{equation} \label{eq318a}
\tilde{T}_d (\pi) | \hat{\mathbf{W}}_0^d = \mathbf{x}^d \convp\pi
\exp(f^\ast l/2) \qquad\mbox{as } \dr.\vspace*{-3pt}
\end{equation}
\end{lem}
\begin{pf}
Fix a sequence $\{ \mathbf{x}^d \}$. Let $\tilde{T}_d^1
(\pi) = \frac{1}{[d^\delta]} \sum_{j=0}^{[d^\beta-1]} \Omega_d
(\hat{\mathbf{W}}_j^d)^{-1}$ and let $\tilde{T}_d^2 (\pi) =
\frac{1}{[d^\delta]} \sum_{j= [d^\beta]}^{[\pi d^\delta-1]}
\Omega_d (\hat{\mathbf{W}}_j^d)^{-1}$. Thus $\tilde{T}_d (\pi) =
\tilde{T}_d^1 (\pi) + \tilde{T}_d^2 (\pi)$.

Let $A_d = \sum_{j=0}^{[d^\delta]} \{ \hat{\mathbf{W}}_j^d \notin
F_d^1 \}$. By Theorem~\ref{lem321}, (\ref{eqss47}), $\pz(A_d |
\hat{\mathbf{W}}_0^d = \mathbf{x}^d) \rightarrow0$ as $\dr$ and
conditional upon $A_d^C$, $\tilde{T}_d^1 (\pi) \leq\frac{[d^\beta]
d^\gamma}{[d^\delta]}$. Hence $\tilde{T}_d^1 (\pi) |
\hat{\mathbf{W}}_0^d = \mathbf{x}^d \convp0$ as $\dr$.

By Lemma~\ref{lem316}, (\ref{eq316a}),
%
%
\begin{eqnarray} \label{eq318b}
\ez[ \tilde{T}_d^2 (\pi) | \hat{\mathbf{W}}_0^d = \mathbf{x}^d
] &=& \frac{1}{[d^\delta]} \sum_{j=[d^\beta]}^{[\pi d^\delta
-1]} \ez[ \Omega_d (\hat{\mathbf{W}}_j^d)^{-1} |
\hat{\mathbf{W}}_0^d = \mathbf{x}^d ] \nonumber\\[-9pt]\\[-9pt]
& \rightarrow& \pi\exp( f^\ast l /2).\nonumber\vadjust{\goodbreak}
\end{eqnarray}
By
Chebyshev's inequality, for any $\varepsilon> 0$,
%
%
\begin{eqnarray} \label{eq318c}\qquad
& & \pz\bigl( \bigl| \tilde{T}_d^2 (\pi) -\ez[ \tilde{T}_d^2
(\pi) | \hat{\mathbf{W}}_0^d = \mathbf{x}^d ] \bigr| >
\varepsilon| \hat{\mathbf{W}}_0^d = \mathbf{x}^d \bigr)
\nonumber\\[-8pt]\\[-8pt]
&&\qquad\leq\frac{1}{\varepsilon^2 [d^\delta]^2} \sum_{j=[d^\beta
]}^{[\pi
d^\delta-1]} \sum_{l=[d^\beta]}^{[\pi d^\delta-1]} \operatorname{cov} \bigl(
\Omega_d (\hat{\mathbf{W}}_j^d)^{-1}, \Omega_d
(\hat{\mathbf{W}}_l^d)^{-1} | \hat{\mathbf{W}}_0^d = \mathbf{x}^d
\bigr).\nonumber
\end{eqnarray}
Since for all $j, l$,
%
%
\begin{eqnarray} \label{eq318d} & & \operatorname{cov} \bigl(
\Omega_d (\hat{\mathbf{W}}_j^d)^{-1}, \Omega_d
(\hat{\mathbf{W}}_l^d)^{-1} | \hat{\mathbf{W}}_0^d = \mathbf{x}^d
\bigr) \nonumber\\[-8pt]\\[-8pt]
&&\qquad\leq \operatorname{var} \bigl( \Omega_d
(\hat{\mathbf{W}}_j^d)^{-1} | \hat{\mathbf{W}}_0^d = \mathbf{x}^d
\bigr)^{{1}/{2}} \operatorname{var} \bigl( \Omega_d
(\hat{\mathbf{W}}_l^d)^{-1} | \hat{\mathbf{W}}_0^d = \mathbf{x}^d
\bigr)^{{1}/{2}},\nonumber
\end{eqnarray}
it is straightforward to show,
using Lemma~\ref{lem317}, that the right-hand side of~(\ref{eq318c})
converges to 0 as $\dr$. Thus $\tilde{T}_d^2 (\pi) |
\hat{\mathbf{W}}_0^d = \mathbf{x}^d \convp\pi\exp(f^\ast l/2)$ as
$\dr$ and the lemma follows immediately.
\end{pf}
%
%
\begin{theorem}
\label{lem319} For any sequence $\{ \mathbf{x}^d \}$ such that
$\mathbf{x}^d \in\tilde{F}_d$,
%
%
\begin{equation} \label{eq319a}
P_d | \mathbf{X}_0^d = \mathbf{x}^d \convp\exp(-f^\ast l/2)
\qquad\mbox{as } \dr.
\end{equation}
\end{theorem}
\begin{pf}
For any $0 < \pi\leq1$, by Lemmas~\ref{lem312},
\ref{lem313} and~\ref{lem318},
%
%
\begin{equation} \label{eq319b} T_d
(\pi) | \hat{\mathbf{X}}_0^d = \mathbf{x}^d \convp\pi\exp(f^\ast
l/2).
\end{equation}

Since $P_d$ satisfies $T_d (P_d) \leq1 < T_d (P_d + 1/[d^\delta])$,
for any $\varepsilon> 0$,
%
%
\begin{eqnarray} \label{eq319c}
& & \pz\bigl(| P_d - \exp(-f^\ast l/2)| > \varepsilon| \mathbf{X}_0^d =
\mathbf{x}^d\bigr) \nonumber\\
&&\qquad\leq\pz\bigl(T_d \bigl(\exp(-f^\ast l/2) -
\varepsilon/2\bigr)
> 1 | \hat{\mathbf{X}}_0^d = \mathbf{x}^d\bigr)\\
&&\qquad\quad{} + \pz\bigl(T_d \bigl(\exp(-f^\ast
l/2) + \varepsilon/2\bigr) \leq1 | \hat{\mathbf{X}}_0^d = \mathbf{x}^d\bigr)
\nonumber
\end{eqnarray}
for all sufficiently large $d$. The lemma follows, since (\ref
{eq319b}) ensures that
the right-hand side of (\ref{eq319c}) converges to 0 as
$\dr$.
\end{pf}

%
%

\section{\texorpdfstring{Proof of (\lowercase{\protect\ref{eqn115}})}{Proof of (2.13)}} \label{SecGen}

From Appendix~\ref{SecPd}, we have that for any sequence
$\{\mathbf{x}^d \}$, such that $\mathbf{x}^d \in\tilde{F}_d$, $P_d
| \mathbf{X}_0^d = \mathbf{x}^d \convp\exp(-lf^\ast/2)$ as $\dr$.
Therefore we proceed by showing that, for any $0 \leq\pi\leq1$,
%
%
\begin{equation} \label{eq41e} \sup_{\mathbf{x}^d \in F_d} | \hat
{G}_d^{\delta, \pi} H
(\mathbf{x}^d) - \pi\hat{G} H (x_1) | \rightarrow0 \qquad
\mbox{as } \dr,
\end{equation}
where $\hat{G}_d^{\delta, \pi} H
(\mathbf{x}^d) = \frac{d^2}{[d^\delta]} \ez[
(H(\hat{\mathbf{X}}^d_{[\pi d^\delta]}) - H (\hat{\mathbf{X}}^d_0))
| \hat{\mathbf{X}}^d_0 = \mathbf{x}^d ]$ is defined in
(\ref{eqmain6}) and
%
%
\begin{equation} \label{eq41a}
\hat{G} H(x) = \frac{l^2}{3} \biggl\{ \frac{1}{2} g'(x) H'(x) +
\frac{1}{2} H''(x) \biggr\}.
\end{equation}
Equation (\ref{eqn115}) will then be proved using the triangle
inequality.\vadjust{\goodbreak}

We analyze $\hat{G}_d H (\xx_j^d) = d^2 \ez[ H (\xx_1^d - \xx_0^d)
| \xx_0^d = \mathbf{x}^d]$, which is defined in (\ref{eqmain6x}),
before using (\ref{eqmain6}) to study $\hat{G}_d^{\delta, \pi} H
(\mathbf{x}^d)$. However,\vspace*{1pt} first we require some definitions and
preliminary results. Throughout we will utilize the following key
facts noted in Section~\ref{SecAlg}: $H^\prime(0) = H^\prime(1) =0$
and that $H^\ast_1, H^\ast_2
<\infty$, where $H^\ast_1 = {\sup_{0 \leq y \leq1}} |H^{\prime} (y)|
$ and $H^\ast_2 = {\sup_{0 \leq y \leq1}} |H^{\prime\prime} (y)|
$.\vspace*{1pt}

We follow~\cite{BR00} and~\cite{NR06} in noting that, for any
function $h$ which is a twice differentiable function on
$\mathbb{R}$, the function $z \mapsto1 \wedge e^{h (z)}$ is also
twice differentiable, except at a countable number of points, with
first derivative given Lebesgue almost everywhere by the function
\[
\frac{d}{dz} 1 \wedge e^{h(z)} = \cases{
h'(z) e^{h(z)}, &\quad if $h(z) < 0$, \cr
0, &\quad if $h(z) \geq0$.}
\]
The second derivative can similarly be
obtained but will not be explicitly required for our calculations.

For $-1 \leq z \leq1$, let $J_d^z (\mathbf{x}^d)$ denote the
probability of accepting a move in the RWM algorithm given that
$Z_{1,1} =z$ and let
%
%
\begin{eqnarray} \label{eqn4a}
\tilde{J}_d^0 (\mathbf{x}^d) &=& \ez
\Biggl[ \exp\Biggl( \sum_{j=2}^d \{ g (x_j + \sigma_d Z_{1,j}) - g
(x_j) \} \Biggr)\nonumber\\[-8pt]\\[-8pt]
&&\hspace*{10.3pt}{}\times 1_{\{\sum_{j=2}^d ( g (x_j + \sigma_d Z_{1,j}) - g
(x_j)) < 0 \}} \prod_{j=2}^d 1_{ \{ 0 < x_j + \sigma_d Z_{1,j} < 1
\}} \Biggr].\nonumber
\end{eqnarray}
Then for all $\mathbf{x}^d$, using Taylor's theorem,
%
%
\begin{equation} \label{eqn4b}
J_d^z (\mathbf{x}^d) = 1_{ \{ 0 < x_1 + \sigma_d z < 1 \}} \{
J_d^0 (\mathbf{x}^d) + \sigma_d g^\prime(x_1) z \tilde{J}_d^0
(\mathbf{x}^d) + O (\sigma_d^2) \}.
\end{equation}
Therefore
for $x_1 \in(\sigma_d, 1- \sigma_d)$,
%
%
\begin{equation} \label{eqn4c} J_d (\mathbf{x}^d) = J_d^0
(\mathbf{x}^d) + O (\sigma_d^2).
\end{equation}

%
\begin{lem}
\label{lem40}
%
%
\begin{equation}
\label{eq40a} \sup_{\mathbf{x}^d \in F_d} \biggl| \frac{\tilde{J}_d^0
(\mathbf{x}^d)}{J_d^0 (\mathbf{x}^d)} - \frac{1}{2} \biggr|
\rightarrow0 \qquad\mbox{as } \dr.
\end{equation}
\end{lem}
\begin{pf}
Let $\tilde{\Omega}_d^0 (\mathbf{x}^d) = \ez[
\prod_{j=2}^d 1_{ \{ 0 < x_j + \sigma_d Z_{1,j} < 1 \}}
1_{\{\sum_{j=2}^d ( g (x_j + \sigma_d Z_{1,j}) - g (x_j)) < 0 \}}
]$ and let $\Omega_d^0 (\mathbf{x}^d) = \ez[
\prod_{j=2}^d 1_{ \{ 0 < x_j + \sigma_d Z_{1,j} < 1 \}} ]$,
the probability\vspace*{1pt} a proposed move stays inside the unit cube given
that the first component does not move. The proof of (\ref{eq33b})
can be adapted to show that, for any $\alpha< \frac{1}{2}$,
$d^\alpha| \Omega_d^0 (\mathbf{x}^d) - J_d^0 (\mathbf{x}^d)|,
d^\alpha| \tilde{\Omega}_d^0 (\mathbf{x}^d) - \tilde{J}_d^0
(\mathbf{x}^d)| \rightarrow0$ as $\dr$. Therefore since for
$\mathbf{x}^d \in F_d$, $J_d^0 (\mathbf{x}^d), \Omega_d^0
(\mathbf{x}^d) \geq\exp(-l g^\ast) d^{-\gamma}$,
(\ref{eqlowerJ}), we have that
%
%
\begin{equation}
\label{eq40b} \sup_{\mathbf{x}^d \in F_d} \biggl| \frac{\tilde{J}_d^0
(\mathbf{x}^d)}{J_d^0 (\mathbf{x}^d)} - \frac{\tilde{\Omega}_d^0
(\mathbf{x}^d)}{\Omega_d^0 (\mathbf{x}^d)} \biggr| \rightarrow0
\qquad\mbox{as } \dr.
\end{equation}
Let $\mathcal{B}_d (\mathbf{x}^d) = \{ 2 \leq j \leq d; x_j \in
R_d^l \}$ and let $I_d (\mathbf{x}^d) = \sum_{j \notin
\mathcal{B}_d (\mathbf{x}^d)} \sigma_d g^\prime(x_j) Z_{1,j}$.
Since $|\mathcal{B}_d (\mathbf{x}^d)| \leq\gamma\log d$, we have
that
%
%
\begin{equation}
\label{eq40ba} \biggl| \sum_{j \in\mathcal{B}_d (\mathbf{x}^d)} \bigl(
g (x_j + \sigma_d Z_{1,j}) - g (x_j)\bigr) \biggr| \leq(\gamma\log d)
\sigma_d g^\ast.
\end{equation}
Then using a Taylor series expansion, there exists $K < \infty$
such that, for all $\mathbf{x}^d \in F_d$,
%
%
\begin{equation}
\label{eq40c}\quad I_d (\mathbf{x}^d) - \frac{K \log d}{d} \leq
\sum_{j=2}^d \bigl( g (x_j + \sigma_d Z_j) - g (x_j)\bigr) \leq I_d
(\mathbf{x}^d) + \frac{K \log d}{d}.
\end{equation}
Since $Z_{1,1},
Z_{1,2}, \ldots,$ are independent, and whether or not a proposed
move from $\mathbf{x}^d$ stays inside the hypercube depends only
upon $\mathcal{B}^d (\mathbf{x}^d)$,
%
%
\begin{eqnarray}
\label{eq40d} &&\Omega_d^0 (\mathbf{x}^d) \pz\bigl(I_d (\mathbf{x}^d) < -
K \log d /d\bigr)\nonumber\\[-8pt]\\[-8pt]
&&\qquad \leq\tilde{\Omega}_d^0 (\mathbf{x}^d) \leq\Omega_d^0
(\mathbf{x}^d) \pz\bigl(I_d (\mathbf{x}^d) < K \log d
/d\bigr).\nonumber
\end{eqnarray}
For all $\mathbf{x}^d \in F_d$, $\frac{1}{d} \sum_{j=1}^d
g^\prime(x_j)^2 \rightarrow\ez[ g^\prime(X_1)^2]$, so
\[
\sqrt{d}
I_d (\mathbf{x}^d) \convd N(0, \ez[ g^\prime(X_1)^2])\qquad \mbox{as $\dr$}.
\]
Therefore it follows that
%
%
\begin{equation}
\label{eq40e} \sup_{\mathbf{x}^d \in F_d } \biggl|
\frac{\tilde{\Omega}_d^0 (\mathbf{x}^d)}{\Omega_d^0 (\mathbf{x}^d)}
- \frac{1}{2} \biggr| \rightarrow0 \qquad\mbox{as } \dr
\end{equation}
with the lemma following from (\ref{eq40b}) and
(\ref{eq40e}) by the triangle inequality.
\end{pf}
%
%
\begin{lem}
\label{lem41} For $x_1 \in(\sigma_d, 1- \sigma_d)$ and
$\mathbf{x}^d \in F_d$,
%
%
\begin{equation}
\label{eq41fx} \hat{G}_dH(\mathbf{x}^d)=
\frac{l^2}{3} \biggl\{ \frac{1}{2} H^{\prime\prime} (x_1) +
\frac{\tilde{J}_d^0 (\mathbf{x}^d)}{J_d^0 (\mathbf{x}^d)} g^\prime
(x_1) H^\prime(x_1) \biggr\} + \varepsilon_d,
\end{equation}
where $\varepsilon_d \rightarrow0$ as $\dr$.

For $x_1 \in R_d^l$,
%
%
\begin{equation} \label{eq41fa} | \hat{G}_d H (\mathbf{x}^d) | \leq
\tfrac{3}{2} H^\ast_2 l^2.
\end{equation}
\end{lem}
\begin{pf}
For $d \geq1$, fix $\mathbf{x}^d \in F_d$ and suppose
that $x_1 \in(\sigma_d, 1-\sigma_d)$. Then
%
%
\begin{eqnarray} \label{eq41ex}\hspace*{32pt}
\hat{G}_d H (\mathbf{x}^d) &=& d^2 \ez[ H (\xx^d_1) - H
(\xx^d_0) | \xx^d_0 = \mathbf{x}^d ] \nonumber\\[-4pt]\\[-12pt]
&=& \frac{d^2}{J_d (\mathbf{x}^d)} \ez\biggl[ \bigl(H (\mathbf{x}^d +
\sigma_d \mathbf{Z}^d) - H (\mathbf{x}^d)\bigr) \biggl\{ 1 \wedge
\frac{\pi_d (\mathbf{x}^d + \sigma_d \mathbf{Z}^d)}{\pi_d
(\mathbf{x}^d)} \biggr\} \biggr].\nonumber
\end{eqnarray}
The right-hand side of (\ref{eq41ex}) is familiar in that it is the
generator of the RWM-algorithm divided by the acceptance
probability; see, for example,~\cite{RGG}, page 113.\vadjust{\goodbreak}

First, note that
\begin{eqnarray*}
H (x_1 + \sigma_d Z_1 ) - H (x_1) &=& \sigma_d Z_1 H' (x_1) +
\frac{\sigma_d^2}{2} Z_1^2 H'' (x_1) \\
&&{}+ \frac{\sigma_d^2}{2} Z_1^2
\{H^{ \prime\prime} (x_1 + \psi_1^d) - H''(x_1) \}.
\end{eqnarray*}
Using
(\ref{eqn4b}), (\ref{eqn4c}) and noting that $0 < x_1 + \sigma_d Z_1
<1$, we have that
%
%
\begin{eqnarray} \label{eq41eb}
\hat{G}_d H (\mathbf{x}^d) & =& \frac{d^2}{J_d^0 (\mathbf{x}^d) + O
(\sigma_d^2)} \nonumber\\
&&\hspace*{0pt}{}\times\ez\biggl[ \biggl\{ \sigma_d Z_1 H' (x_1) +
\frac{\sigma_d^2}{2} Z_1^2 H'' (x_1)\nonumber\\
&&\hspace*{29pt}{} + \frac{\sigma_d^2}{2} Z_1^2
\{H''(x_1 + \psi_1^d) - H''(x_1) \} \biggr\} \nonumber\\
&&\hspace*{24.2pt}{} \times\{ J_d^0 (\mathbf{x}^d) +
\tilde{J}_d^0 (\mathbf{x}^d) \sigma_d g'(x_1) Z_1 + O (\sigma_d^2)
\} 1_{\{0 < x_1 + \sigma_d Z_1 <1 \}} \biggr] \nonumber\\
&=& \frac{d^2 J_d^0 (\mathbf{x}^d)}{J_d^0 (\mathbf{x}^d) + O
(\sigma_d^2)} \sigma_d \ez[Z_1] H^\prime(x_1)\nonumber\\
&&{} + \frac{d^2 J_d^0
(\mathbf{x}^d)}{J_d^0 (\mathbf{x}^d) + O (\sigma_d^2)}
\frac{\sigma_d^2}{2} \ez[Z_1^2] H^{ \prime\prime} (x_1)
\nonumber\\
&&{}+ \frac{d^2 J_d^0 (\mathbf{x}^d)}{J_d^0 (\mathbf{x}^d) + O
(\sigma_d^2)} \frac{\sigma_d^2}{2} \ez[Z_1^2 \{H^{ \prime\prime}
(x_1 + \psi_1^d) - H''(x_1) \} ] \\
&&{} + \frac{d^2 \tilde{J}_d^0 (\mathbf{x}^d)}{J_d^0 (\mathbf{x}^d) +
O (\sigma_d^2)} \sigma_d^2 g^\prime(x_1) H^\prime(x_1) \ez[Z_1^2]\nonumber\\
&&{}+ \frac{d^2 }{J_d^0 (\mathbf{x}^d) + O (\sigma_d^2)} O
(\sigma_d^3). \nonumber
\end{eqnarray}
The first term on the right-hand side of (\ref{eq41eb}) is 0. Since
$H^\ast_2 < \infty$, by the continuous mapping theorem, $\{H^{
\prime\prime} (x_1 + \psi_1^d) - H''(x_1) \} \convp0$ as $\dr$ and
then since $Z_1$ is bounded the third term on the right-hand side of
(\ref{eq41eb}) converges to 0 as $\dr$. For $\mathbf{x}^d \in F_d$,
$J_d^0 (\mathbf{x}^d) \geq e^{-lg^\ast} d^{- \gamma}$, and so, the
right-hand side of (\ref{eq41eb}) equals
\[
\frac{l^2}{3} \biggl\{ \frac{1}{2} H^{\prime\prime} (x_1) +
\frac{\tilde{J}_d^0 (\mathbf{x}^d)}{J_d^0 (\mathbf{x}^d)} g^\prime
(x_1) H^\prime(x_1) \biggr\} + \varepsilon_d,
\]
where $\varepsilon_d \rightarrow0$ as $\dr$. Thus (\ref{eq41fx})
is proved.

The proof of (\ref{eq41fa}) follows straightforwardly using Taylor
series expansions since $H^\prime(0) = H^\prime(1)=0$.\vadjust{\goodbreak}
\end{pf}

Since $g^\ast= {\sup_{0 \leq y \leq1}} |g^\prime(y)|, H^\ast_1,
H^\ast_2 < \infty$, an immediate
consequence of Lemma~\ref{lem41} is that, there exists $K^\ast<
\infty$ such that
%
%
\begin{equation} \label{eqn4d} {\sup_d \sup_{\mathbf{x}^d \in F_d}}
|\hat{G}_d H
(\mathbf{x}^d)| \leq K^\ast.\vspace*{-2pt}
\end{equation}

%
\begin{lem}
\label{lem42} For any sequence of positive integers $\{ k_d \}$
satisfying $[d^\beta] \leq k_d \leq[d^\delta]$,
%
%
\begin{equation}\label{eq42aa}
\sup_{\mathbf{x}^d \in\tilde{F}_d} \bigl| \ez[ \hat{G}_d H
(\hat{\mathbf{X}}_{k_d}^d) | \hat{\mathbf{X}}_0^d = \mathbf{x}^d] -
\hat{G} H (x_1) \bigr| \rightarrow0 \qquad\mbox{as } \dr.\vspace*{-2pt}
\end{equation}
\end{lem}
\begin{pf}
Fix $\{k_d \}$ and note that
%
%
\begin{eqnarray} \label{eq42ab}
&&
\ez[ \hat{G}_d H (\xx_{k_d}^d) | \hat{\mathbf{X}}_0^d =
\mathbf{x}^d ] \nonumber\\
&&\qquad= \pz(\xx_{k_d}^d \in F_d |
\hat{\mathbf{X}}_0^d = \mathbf{x}^d) \ez[ \hat{G}_d H
(\xx_{k_d}^d) | \hat{\mathbf{X}}_0^d =
\mathbf{x}^d, \xx_{k_d}^d \in F_d ] \\
&&\qquad\quad{} + \pz(\xx_{k_d}^d \notin F_d | \hat{\mathbf{X}}_0^d =
\mathbf{x}^d) \ez[ \hat{G}_d H (\xx_{k_d}^d) |
\hat{\mathbf{X}}_0^d = \mathbf{x}^d, \xx_{k_d}^d \notin F_d
].\nonumber
\end{eqnarray}
Since $H \in\mathcal{D}$, $H^\ast_0 = \sup_{0 \leq y \leq1} |
H(y)| < \infty$. Therefore, for all $\mathbf{y}^d \in[0,1]^d$,
$\hat{G}_d H (\mathbf{y}^d) \leq2 d^2 H^\ast_0$. By
(\ref{emainx1}), $\sup_{\mathbf{x}^d \in\tilde{F}_d} d^2 \pz
(\hat{\mathbf{X}}_{k_d}^d \notin F_d | \xx_0^d = \mathbf{x}^d)
\rightarrow0$ as $\dr$. Thus the latter term on the right-hand side
of (\ref{eq42ab}) converges to 0 as $\dr$.

Now
%
%
\begin{eqnarray} \label{eq42ac} & & \ez[ \hat{G}_d H
(\xx_{k_d}^d) | \hat{\mathbf{X}}_0^d = \mathbf{x}^d, \xx_{k_d}^d
\in
F_d ] \nonumber\\
&&\qquad= \pz(\hat{X}_{k_d,1}^d \notin R_d^l | \hat{\mathbf{X}}_0^d =
\mathbf{x}^d, \xx_{k_d}^d \in F_d) \nonumber\\
&&\qquad\quad{}\times\ez[ \hat{G}_d H
(\xx_{k_d}^d) | \hat{\mathbf{X}}_0^d = \mathbf{x}^d, \xx_{k_d}^d
\in
F_d, \hat{X}_{k_d,1}^d \notin R_d^l ] \\
&& \qquad\quad{} + \pz(\hat{X}_{k_d,1}^d \in R_d^l | \hat{\mathbf
{X}}_0^d =
\mathbf{x}^d, \xx_{k_d}^d \in F_d) \nonumber\\
&&\qquad\quad\hspace*{11pt}{}\times\ez[ \hat{G}_d H
(\xx_{k_d}^d) | \hat{\mathbf{X}}_0^d = \mathbf{x}^d, \xx_{k_d}^d
\in
F_d, \hat{X}_{k_d,1}^d \in R_d^l ].\nonumber
\end{eqnarray}
Consider
first the latter term on the right-hand side of (\ref{eq42ac}). By
Lemma~\ref{lem41}, (\ref{eq41fa}),
%
%
\begin{equation} \label{eq42aca}
\ez[ \hat{G}_d H (\xx_{k_d}^d) | \hat{\mathbf{X}}_0^d =
\mathbf{x}^d, \xx_{k_d}^d \in F_d, \hat{X}_{k_d,1}^d \in R_d^l
] \leq\tfrac{3}{2} l^2 H^\ast_2.
\end{equation}
Note that
%
%
\begin{eqnarray} \label{eq42ad}\qquad
\pz(\hat{X}_{k_d,1}^d \in R_d^l | \hat{\mathbf{X}}_0^d =
\mathbf{x}^d, \xx_{k_d}^d \in F_d) &=& \frac{\pz(\hat{X}_{k_d,1}^d
\in R_d^l, \xx_{k_d}^d \in F_d | \hat{\mathbf{X}}_0^d =
\mathbf{x}^d)}{\pz(\xx_{k_d}^d \in F_d |\hat{\mathbf{X}}_0^d =
\mathbf{x}^d)} \nonumber\\[-8pt]\\[-8pt]
& \leq& \frac{\pz(\hat{X}_{k_d,1}^d
\in R_d^l | \hat{\mathbf{X}}_0^d = \mathbf{x}^d)}{\pz(\xx_{k_d}^d
\in F_d |\hat{\mathbf{X}}_0^d = \mathbf{x}^d)}.\nonumber
\end{eqnarray}

By (\ref{emainx1}), for $\mathbf{x}^d \in\tilde{F}_d$, $\pz
(\xx_{k_d}^d \in F_d |\hat{\mathbf{X}}_0^d = \mathbf{x}^d)
\rightarrow1$ as $\dr$. Use Corollary~\ref{lem35} and Lemma\vadjust{\goodbreak}
\ref{lem38a} to show that $\pz(\hat{X}_{k_d,1}^d \in R_d^l |
\hat{\mathbf{X}}_0^d = \mathbf{x}^d) \rightarrow0$ as $\dr$. Hence,
the right-hand side of (\ref{eq42ad}) converges to 0 as $\dr$ and
consequently the latter term on the right-hand side of (\ref{eq42ac})
converges to 0 as $\dr$.

It follows from the above arguments that
%
%
\begin{equation}
\label{eq42ae} \min_{\mathbf{x}^d \in\tilde{F}_d} \pz
(\hat{X}_{k_d,1}^d \notin R_d^l, \xx_{k_d}^d \in F_d |
\hat{\mathbf{X}}_0^d = \mathbf{x}^d) \rightarrow1 \qquad
\mbox{as } \dr.
\end{equation}
%
%
Also it follows from (\ref{eqn4d}) that there exists $K < \infty$
such that
%
%
\begin{equation} \label{eq42af}
\sup_d \sup_{\mathbf{x}^d \in\tilde{F}_d} \ez[ \hat{G}_d H
(\xx_{k_d}^d) | \hat{\mathbf{X}}_0^d = \mathbf{x}^d, \xx_{k_d}^d
\in
F_d, \hat{X}_{k_d,1}^d \notin R_d^l ] \leq K.
\end{equation}
Therefore, it is straightforward using (\ref{eq42ab}),
(\ref{eq42ac}) and the triangle inequality to show that
%
%
\begin{eqnarray} \label{eq42ag} && \sup_{\mathbf{x}^d \in\tilde{F}_d}
\bigl| \ez[ \hat{G}_d H (\xx_{k_d}^d) | \hat{\mathbf
{X}}_0^d =
\mathbf{x}^d ] \nonumber\\
&&\qquad\quad\hspace*{-8.6pt}{}- \ez[ \hat{G}_d H (\xx_{k_d}^d) |
\hat{\mathbf{X}}_0^d = \mathbf{x}^d, \xx_{k_d}^d \in F_d,
\hat{X}_{k_d,1}^d \notin R_d^l ] \bigr|\\
&& \qquad \rightarrow0
\qquad\mbox{as } \dr.\nonumber
\end{eqnarray}

By Lemma~\ref{lem41}, (\ref{eq41fx}), there exists $\varepsilon_d^1
\rightarrow0$ as $\dr$, such that
%
%
\begin{eqnarray} \label{eq42ah} & & \sup_{\mathbf{x}^d \in\tilde{F}_d}
\bigl| \ez[ \hat{G}_d H (\xx_{k_d}^d) - \hat{G} H
(\hat{X}_{k_d,1}^d) | \hat{\mathbf{X}}_0^d = \mathbf{x}^d,
\xx_{k_d}^d \in F_d, \hat{X}_{k_d,1}^d \notin R_d^l ] \bigr|
\nonumber\\[-1pt]
&&\qquad\leq\frac{l^2}{3} \sup_{0 \leq y \leq1} | g^\prime(y)
H^\prime(y)| \nonumber\\[-8.5pt]\\[-8.5pt]
&&\qquad\quad{}\times\sup_{\mathbf{x}^d \in\tilde{F}_d} \ez\biggl[ \biggl|
\frac{\tilde{J}_d^0 (\xx_{k_d}^d)}{J_d^0 (\xx_{k_d}^d)} -
\frac{1}{2} \biggr| \Big| \hat{\mathbf{X}}_0^d = \mathbf{x}^d,
\xx_{k_d}^d \in F_d, \hat{X}_{k_d,1}^d \notin R_d^l \biggr] +
\varepsilon_d^1 \nonumber\\[-1pt]
&&\qquad\leq\frac{l^2}{3} g^\ast H^\ast_1 \sup_{\mathbf{y}^d \in F_d}
\biggl| \frac{\tilde{J}_d^0 (\mathbf{y}^d)}{J_d^0 (\mathbf{y}^d)} -
\frac{1}{2} \biggr| + \varepsilon_d^1.\nonumber
\end{eqnarray}
By Lemma~\ref{lem40}, the right-hand side of (\ref{eq42ah})
converges to 0 as $\dr$.

Using the triangle inequality, the lemma follows by showing that
%
%
\begin{eqnarray} \label{eq42ai}
&&\sup_{\mathbf{x}^d \in\tilde{F}_d}
\bigl| \ez[ \hat{G} H (\hat{X}_{k_d,1}^d) |
\hat{\mathbf{X}}_0^d = \mathbf{x}^d, \xx_{k_d}^d \in F_d,
\hat{X}_{k_d,1}^d \notin R_d^l ] - \hat{G} H (x_1)
\bigr|\nonumber\\[-8.5pt]\\[-8.5pt]
&&\qquad\rightarrow0 \qquad\mbox{as } \dr.\nonumber
\end{eqnarray}
Note
that $|\hat{X}_{k_d,1}^d -x_1| \leq k_d \sigma_d$, and so,
(\ref{eq42ai}) follows since $\hat{G} H (\cdot)$ is continuous.
\end{pf}

We are in position to prove (\ref{eq41e}).
%
%
\begin{lem}
\label{lem43} For any $0 \leq\pi\leq1$,
%
%
\begin{equation}\label{eq43a}
\sup_{\mathbf{x}^d \in\tilde{F}_d} | \hat{G}_d^{\delta, \pi}
(\mathbf{x}^d) - \pi\hat{G} H (x_1) | \rightarrow0 \qquad
\mbox{as } \dr.\vadjust{\goodbreak}
\end{equation}
\end{lem}
\begin{pf}
Since (\ref{eq43a}) trivially holds for $\pi=0$, we
assume that $\pi> 0$. For all sufficiently large $d$, by the
triangle inequality,
%
%
\begin{eqnarray} \label{eq43b} & & | \hat{G}_d^{\delta, \pi}
(\mathbf{x}^d) - \pi\hat{G} H
(x_1) | \nonumber\\[-1pt]
&&\qquad= \Biggl| \frac{1}{[d^\delta]} \sum_{j=0}^{[\pi d^\delta-1]}
\ez
[ \hat{G}_d H (\xx_j^d) | \hat{\mathbf{X}}_0^d = \mathbf{x}^d
] - \pi\hat{G} H (x_1) \Biggr| \nonumber\\[-1pt]
&&\qquad\leq\Biggl| \frac{1}{[d^\delta]} \sum_{j=0}^{[d^\beta] -1} \ez
[ \hat{G}_d H (\hat{\mathbf{X}}_j^d) | \hat{\mathbf{X}}_0^d =
\mathbf{x}^d ] \Biggr| \\[-1pt]
& &\qquad\quad{} + \frac{1}{[d^\delta]} \sum_{j=[d^\beta]}^{[\pi
d^\delta-1]}
\bigl| \ez[ \hat{G}_d H (\hat{\mathbf{X}}_j^d) |
\hat{\mathbf{X}}_0^d = \mathbf{x}^d ] - \hat{G} H (x_1)
\bigr|\nonumber\\[-1pt]
&&\qquad\quad{} + \biggl( \pi- \frac{[\pi d^\delta] - [d^\beta]}{[d^\delta]}
\biggr) \hat{G} H (x_1).\nonumber
\end{eqnarray}

Since
%
%
\begin{eqnarray} \label{eq43c}
&&\ez[ \hat{G}_d H (\hat{\mathbf{X}}_j^d) | \hat{\mathbf{X}}_0^d
= \mathbf{x}^d ] \nonumber\\
&&\qquad =
\ez
[ \hat{G}_d H (\hat{\mathbf{X}}_j^d) | \hat{\mathbf{X}}_0^d =
\mathbf{x}^d , \xx_j^d \in F_d ] \pz(\xx_j^d \in F_d |
\hat{\mathbf{X}}_0^d = \mathbf{x}^d) \\
&&\qquad\quad{} + \ez[ \hat{G}_d H (\hat{\mathbf{X}}_j^d) |
\hat{\mathbf{X}}_0^d = \mathbf{x}^d , \xx_j^d \notin F_d ]
\pz(\xx_j^d \notin F_d | \hat{\mathbf{X}}_0^d = \mathbf{x}^d),\nonumber
\end{eqnarray}
it is straightforward, following a similar argument to the proof of
Lemma~\ref{lem42}, (\ref{eq42af}), to show that there exists
$\tilde{K} < \infty$ such that, for all $0 \leq j \leq[d^\delta]$,
%
%
\begin{equation} \label{eq43d} \sup_{\mathbf{x}^d \in\tilde{F}_d}
\bigl|
\ez[ \hat{G}_d H (\hat{\mathbf{X}}_j^d) | \hat{\mathbf{X}}_0^d
= \mathbf{x}^d ] \bigr| \leq\tilde{K}.
\end{equation}
Therefore the first term on the right-hand side of (\ref{eq43c}) is
bounded by $[d^\beta] \tilde{K}/[d^\delta]$. By Lemma~\ref{lem42}
the supremum over $\mathbf{x}^d \in\tilde{F}_d$ of the second term
on the right-hand side of (\ref{eq43b}) converges to 0 as $\dr$ and
the lemma follows.
\end{pf}
%
%
\begin{coro}
\label{lem44}
%
%
\begin{equation} \label{eq44a}
\sup_{0 \leq\pi\leq1} \sup_{\mathbf{x}^d \in\tilde{F}_d} |
\hat{G}_d^{\delta, \pi} (\mathbf{x}^d) - \pi\hat{G} H (x_1) |
\rightarrow0 \qquad\mbox{as } \dr.
\end{equation}
\end{coro}
\begin{pf}
Fix $\varepsilon> 0$ and let $\Pi_\varepsilon= \{ 0,
\varepsilon, 2 \varepsilon, \ldots, [1/\varepsilon] \varepsilon, 1\}$. It
follows from Lemma~\ref{lem43} that, for all sufficiently large $d$,
%
%
\begin{equation} \label{eq44b}
\max_{\pi\in\Pi_\varepsilon} \sup_{\mathbf{x}^d \in\tilde{F}_d} |
\hat{G}_d^{\delta, \pi} (\mathbf{x}^d) - \pi\hat{G} H (x_1) |
\leq
\varepsilon.
\end{equation}

Consider any $0 \leq\pi\leq1$. There exists $\tilde{\pi} \in
\Pi_\varepsilon$ such that $\tilde{\pi} \leq\pi< \tilde{\pi} +
\varepsilon$. By the triangle inequality,
%
%
\begin{eqnarray} \label{eq44c}
& &| \hat{G}_d^{\delta, \pi} H (\mathbf{x}^d) - \pi\hat{G} H (x_1)
| \nonumber\\
&&\qquad\leq| \hat{G}_d^{\delta, \pi} H (\mathbf{x}^d) -
\hat{G}_d^{\delta, \tilde{\pi}} H (\mathbf{x}^d)| + |
\hat{G}_d^{\delta, \tilde{\pi}} H (\mathbf{x}^d) - \tilde{\pi}
\hat{G} H (x_1) |\\
&&\qquad\quad{} + (\pi- \tilde{\pi}) |\hat{G} H (x_1)|.\nonumber
\end{eqnarray}
Again by the triangle inequality,
%
%
\begin{eqnarray} \label{eq44d}
&&
\sup_{\mathbf{x}^d \in\tilde{F}_d} | \hat{G}_d^{\delta, \pi
} H
(\mathbf{x}^d) - \hat{G}_d^{\delta, \tilde{\pi}} H (\mathbf{x}^d)
| \nonumber\\[-8pt]\\[-8pt]
&&\qquad\leq\frac{1}{[d^\delta]} \sum_{j = [\tilde{\pi}
d^\delta]}^{[\pi d^\delta-1]} \sup_{\mathbf{x}^d \in\tilde{F}_d}
\bigl|\ez[ \hat{G}_d H (\hat{\mathbf{X}}_j^d) | \hat{\mathbf{X}}_0^d =
\mathbf{x}^d ]\bigr|.\nonumber
\end{eqnarray}
Since for all sufficiently large
$d$, $([\pi d^\delta-1] - [\tilde{\pi} d^\delta])/[d^\delta] \leq2
\varepsilon$, it follows from~(\ref{eq43d}) that the right-hand side of
(\ref{eq44d}) is bounded by $2 \tilde{K} \varepsilon$, where
$\tilde{K}$ is defined in Lemma~\ref{lem43}.

Let $\hat{K} = 2 \tilde{K} +1 + \sup_{0 \leq y \leq1} |\hat{G} H
(y)|$. Note that since $g^\ast, H^\ast_1, H^\ast_2 < \infty$, we
have that
$\hat{K} < \infty$. Therefore it follows from (\ref{eq44c}) that for
all sufficiently large $d$,
%
%
\begin{equation} \label{eq44e} \sup_{\mathbf{x}^d \in\tilde{F}_d}
| \hat{G}_d^{\delta, \pi} H (\mathbf{x}^d) - \pi\hat{G} H (x_1)
| \leq\hat{K} \varepsilon.
\end{equation}
Since (\ref{eq44e}) holds
for all $0 \leq\pi\leq1$ and $\varepsilon> 0$, the lemma follows.
\end{pf}

Finally we are in position to prove (\ref{eqn115}), and hence
complete the proof of Theorem~\ref{main}.
%
%
\begin{lem}
\label{lem45}
%
%
\begin{equation}\label{eq45a}
\sup_{\mathbf{x}^d \in\tilde{F}_d} | G_d^\delta H(\mathbf{x}^d) -
G H (x_1) | \rightarrow0 \qquad\mbox{as } \dr.
\end{equation}
\end{lem}
\begin{pf}
Note that $G_d^\delta H (\mathbf{x}^d)$ is given by
(\ref{eqn113}) and $G H (x_1) = \exp(-l f^\ast/2) \times\hat{G} H
(x_1)$. Therefore by the triangle inequality,
%
%
\begin{eqnarray} \label{eq45ba}\quad & & \sup_{\mathbf{x}^d \in\tilde
{F}_d} | G_d^\delta H(\mathbf{x}^d) - G
H (x_1) | \nonumber\\
&&\qquad= \sup_{\mathbf{x}^d \in\tilde{F}_d} \biggl|
\frac{d^2}{[d^\delta]} \ez\bigl[ H \bigl(\hat{\mathbf{X}}_{[P_d
d^\delta]}^d\bigr) - H (\hat{\mathbf{X}}_0^d) | \hat{\mathbf{X}}_0^d =
\mathbf{x}^d \bigr] - \exp(- l f^\ast/2) \hat{G} H (x_1) \biggr|
\nonumber\\
&&\qquad\leq\sup_{\mathbf{x}^d \in\tilde{F}_d} \biggl|
\ez\biggl[ \frac{d^2}{[d^\delta]} \bigl( H
\bigl(\hat{\mathbf{X}}_{[P_d d^\delta]}^d\bigr) - H (\hat{\mathbf{X}}_0^d)
\bigr) - P_d \hat{G} H (x_1) | \hat{\mathbf{X}}_0^d = \mathbf{x}^d
\biggr] \biggr| \nonumber\\
&&\qquad\quad{}
+ \sup_{\mathbf{x}^d
\in\tilde{F}_d} \bigl| \ez[ P_d \hat{G} H (x_1) |
\hat{\mathbf{X}}_0^d = \mathbf{x}^d ] - \exp
(- l f^\ast/2) \hat{G} H (x_1) \bigr| \nonumber\\[-8pt]\\[-8pt]
&&\qquad\leq\sup_{0 \leq\pi\leq1} \sup_{\mathbf{x}^d \in
\tilde{F}_d} \biggl|
\ez\biggl[ \frac{d^2}{[d^\delta]} \bigl( H
\bigl(\hat{\mathbf{X}}_{[\pi d^\delta]}^d\bigr) - H (\hat{\mathbf{X}}_0^d)
\bigr) - \pi\hat{G} H (x_1) | \hat{\mathbf{X}}_0^d = \mathbf{x}^d
\biggr] \biggr| \nonumber\\
&&\qquad\quad{}+ \sup_{\mathbf{x}^d
\in\tilde{F}_d} \bigl| \ez[ P_d | \hat{\mathbf{X}}_0^d =
\mathbf{x}^d ] - \exp(- l f^\ast/2)
\bigr| \sup_{0 \leq y \leq1} |\hat{G} H (y)| \nonumber\\
&&\qquad\leq\sup_{0 \leq\pi\leq1} \sup_{\mathbf{x}^d \in
\tilde{F}_d} | \hat{G}_d^{\delta, \pi} H (\mathbf{x}^d) -
\pi
\hat{G} H (x_1) | \nonumber\\
&&\qquad\quad{}
+ \sup
_{\mathbf{x}^d \in\tilde{F}_d} \bigl| \ez[ P_d |
\hat{\mathbf{X}}_0^d = \mathbf{x}^d ] - \exp(- l f^\ast/2) \bigr|
\sup_{0 \leq y \leq1} |\hat{G} H (y)|.\nonumber
\end{eqnarray}
By Corollary
\ref{lem44}, the first term on the right-hand side of (\ref{eq45ba})
converges to 0 as $\dr$. By Theorem~\ref{lem319}, for any sequence
$\{ \mathbf{x}^d \}$ such that $\mathbf{x}^d \in\tilde{F}_d$, $P_d
| \mathbf{X}_0^d = \mathbf{x}^d \convp\exp(- l f^\ast/2)$ as
$\dr$. Hence the latter term on the right-hand side of (\ref{eq45ba})
converges to 0 as $\dr$, since $g^\ast, H^\ast_1, H^\ast_2 < \infty$
implies that $\sup_{0 \leq y \leq1} |\hat{G} H (y)| <\infty$.
\end{pf}
\end{appendix}

\section*{Acknowledgments}

We thank the anonymous referees for their helpful comments which have
improved the presentation of the paper.


%

\printaddresses

\end{document}